\documentclass{ws-ijbc}

\usepackage{ws-rotating}     % used only when sideways tables/figures are used
\usepackage{xcolor}
\usepackage[verbose]{hyperref}
\usepackage{tikz}
\usepackage{tikz-cd}
\hypersetup{colorlinks=false,allbordercolors=blue,pdfborderstyle={/S/U/W 1}}
\usetikzlibrary{cd,arrows,babel}
\newcommand\encircle[1]{%
  \tikz[baseline=(X.base)] 
    \node (X) [draw, shape=circle, inner sep=0] {\strut #1};}
\usepackage{mathtools}

\begin{document}

\catchline{}{}{}{}{} % Publisher's Area please ignore

\markboth{N.V. Barabash et al.}{Analytically tractable reconstruction of singular hyperbolic and quasi-strange attractors of Lorenz-type systems}

\title{Analytically tractable reconstruction of singular hyperbolic and quasi-strange attractors of Lorenz-type systems}

\author{Nikita\,V. Barabash}%

\address{Department of Control Theory and Dynamics of Systems, Lobachevsky State University of Nizhny Novgorod,\\
23, Gagarin Ave., 603022, Nizhny Novgorod, Russia\\
Department of Mathematics, Volga State University of Water Transport \\
5, Nesterov Str., 603950, Nizhny Novgorod, Russia\\
barabash@itmm.unn.ru}

\author{Daria\,A.Bakalina}
\address{Department of Control Theory and Dynamics of Systems, Lobachevsky State University of Nizhny Novgorod,\\
23, Gagarin Ave., 603022, Nizhny Novgorod, Russia\\
darya.bakalina@gmail.com}

\author{Vladimir\,N.Belykh}%

\address{Department of Control Theory and Dynamics of Systems, Lobachevsky State University of Nizhny Novgorod,\\
23, Gagarin Ave., 603022, Nizhny Novgorod, Russia\\
Department of Mathematics, Volga State University of Water Transport \\
5, Nesterov Str., 603950, Nizhny Novgorod, Russia\\
belykh@unn.ru}

\maketitle
%\begin{history}
%\received{(May 07, 2025)}; \accepted{(Month Day Year)}; \published{(Month Day Year)}
%\end{history}

\begin{abstract}
We consider a certain three-dimensional piecewise linear system of Lorenz type in the cases of positive and negative saddle value, which is the sum of two eigenvalues of the saddle nearest to zero.
This system was recently proposed and studied in \cite{belykh2019lorenz} for the case of the positive saddle value.
Here we consider the main bifurcations leading to the birth of the strange attractor and provide its comparison with those of original smooth Lorenz model.
For the case of a negative saddle value, we obtain an explicit forms of homoclinic and pitchfork bifurcations forming a codimension 1 cascade leading to the emergence of a quasi-strange attractor.
This cascade reproduces a typical bifurcation route towards the birth of a quasi-strange attractor in smooth Lorenz-like systems with the negative saddle value. We support our analytical result with a numerical comparison with the Lorenz-Lyubimov-Zaks system.
\end{abstract}


%\begin{multicols}{2}
\vspace*{-0.5in}
\section{Introduction}
More than half a century ago, Edward Lorenz discovered complex nonperiodic motions in a system of three deterministic ordinary differential equations, which later received his name \cite{lorenz1963deterministic}.
Since then, the Lorenz system with its strange attractor had become an icon of chaos theory, to which thousands of papers in the field of pure and applied mathematics had been devoted.

The most famous mathematical results include geometric models of the Lorenz attractor \cite{guckenheimer1976strange,afraimovich1977origin,afraimovich1982attractive,williams1979structure}, numerically obtained scenarios of its birth \cite{sparrow1982lorenz}, proofs of the existence of a homoclinic butterfly in Lorenz-type systems \cite{petrovskaya1980homoclinic,belykh1984bifurcation}, asymptotic criteria for the birth of a strange attractor \cite{shilnikov1980appendix,ovsyannikov2016analytic}, models for its destruction \cite{bykov1989boundaries,barrio2012kneadings,creaser2017finding}, as well as a numerical proof of the existence of a strange attractor for canonical values of the parameters from Lorenz's original paper \cite{tucker1999lorenz}.

Traditionally, smooth Lorenz-type flows are associated with the one-dimensional model map 
\begin{equation*}
    \bar{x}=(-\mu+|x|^{\nu})\textrm{\,sign }x,\quad x\in\mathbb{R},\quad\mu\in\mathbb{R},\quad \nu>0,
\end{equation*}
obtained as a result of factorization of the mapping of a 2-D Poincar\'{e} cross-section under the assumption of the existence of a stable foliation on it.
Here, $\mu$ is the splitting value for separatrix loop, and $\nu$ is the saddle index being the absolute value of the ratio of the positive eigenvalue of the saddle to the negative one nearest to zero. 

It is commonly emphasized that this factor map is related to the flow only for small $|\mu|$ and $\nu\approx 1$, i.e. in a small neighborhood of the homoclinic butterfly to a neutral saddle.  
Therefore, despite its long history, the global partitioning of the parameter space of Lorenz-type systems has so far been obtained only in numerical experiments.

In contrast, in this paper, instead of making assumptions or postulating properties that would allow 3-D flows to be confined to a geometric model, we consider a specific 3-D piecewise-linear (PWL) Lorenz-type flow proposed in the paper~\cite{belykh2019lorenz}. 
This system allowed the authors to obtain a two-dimensional Poincar\'{e} map in the form of a skew product and then, using its 1-D Lorenz  factor map, to prove global properties of the flow that were previously inaccessible.

In this paper, we compare these properties with those for the classical Lorenz system, following the main bifurcations and the birth of strange attractors in both systems.

The main result of this paper is the extending of the parameter space into a region with a ``contracting saddle'' achieved when $\nu>1$, i.e. the sum of its closest to zero eigenvalues  becomes negative.

For smooth Lorenz-type systems, this case was considered in a number of papers \cite{arneodo1981possible,lyubimov1983two,rovella1993dynamics,morales2006contracting,kazakov2021bifurcations,araujo2021statistical}, where the birth of a quasi-strange attractor (i.e., containing stable orbits of large period and truly chaotic only on a Cantor set of parameters) has been shown.

However, the obtained results are either numerical or pertain to a small neighborhood of the bifurcation of codimension two, corresponding to a homoclinic butterfly to the neutral saddle.

In this paper, using a one-dimensional factor map for a piecewise linear Lorenz-type flow with a negative saddle value, we obtain global partition of its parameter space.

We show that the obtained analytical result qualitatively reproduces bifurcations in the Lyubimov-Zaks system, which generalizes the Lorenz model to the case of a negative saddle value.

The paper is organized as follows. 
In Section~\ref{sec:model}, we briefly recall a piecewise-linear Lorenz-type system, introduced in ~\cite{belykh2019lorenz}, and corresponding 1-D factor map, derived from it. 
In Section~\ref{sec:positive}, we consider this system in the case of positive saddle value, when the system has codimension 1 bifurcation route leading to birth of Lorenz-type singular hyperbolic strange  attractor. 
We compare it with well-known corresponding bifurcation route of original Lorenz model.  
In Section~\ref{sec:negative}, we consider the case of negative saddle value, for which we prove the bifurcation cascade for Lorenz-type quasi-strange attractor emergence.
We show that this cascade qualitatively coincide with the formation of the quasi-strange attractor in the Lorenz-Lyubimov-Zaks system.

\section{Piecewise-linear Lorenz-type system \label{sec:model}}
Consider a piecewise-linear (PWL) system of Lorenz-type introduced in \cite{belykh2019lorenz}.
The system is composed of three linear subsystems $A_s$, $A_l$ and $A_r$ 
\begin{equation}
A_s :\left\{
\begin{array}{l}
   \dot{x} = x, \\
   \dot{y}  = -\alpha y,\\
   \dot{z}  = -\nu z,
\end{array}\right.
\;
A_{l} : \left\{ 
\begin{array}{l}
   \dot{x} = - \lambda (x+1) + \omega (z - b), 
   \\
   \dot{y}  = - \delta (y + 1),
   \\
   \dot{z}  = -\omega (x + 1) - \lambda (z - b),
\end{array}\right.
\;
A_{r} : \left\{ 
\begin{array}{l}
   \dot{x} = - \lambda (x-1) + \omega (z - b), 
   \\
   \dot{y}  = - \delta (y - 1),
   \\
   \dot{z}  = -\omega (x - 1) - \lambda (z - b),
\end{array}\right.
\label{sys:PWL_Lorenz}
\end{equation}
defined in the following partition of 3-D phase space 
\begin{equation*}
\begin{array}{l}
 G_s=\lbrace|x| < 1, y  \in \mathbb{R}^1, z< b\rbrace,\\
G_l=\lbrace x \leq -1 , (y, z)\in \mathbb{R}^2\rbrace\cup
\lbrace x < 1, y < 0, z > b\rbrace,\\
G_r=\lbrace x \geq 1 , (y, z)\in \mathbb{R}^2\rbrace\cup
\lbrace x > -1, y > 0, z > b\rbrace,
\end{array} 
\end{equation*}
respectively.
All parameters $\alpha$, $\delta$, $\nu$, $\omega$, $\lambda$ and $b$ are positive, and the saddle index $\nu$, explicitly included in the system as a parameter, satisfies the conditions
\begin{equation*}
\frac{1}{2} < \nu < \alpha.
\label{eq:nu<alpha}
\end{equation*}

\begin{figure}
    \centering
    \includegraphics[width=0.32\textwidth]{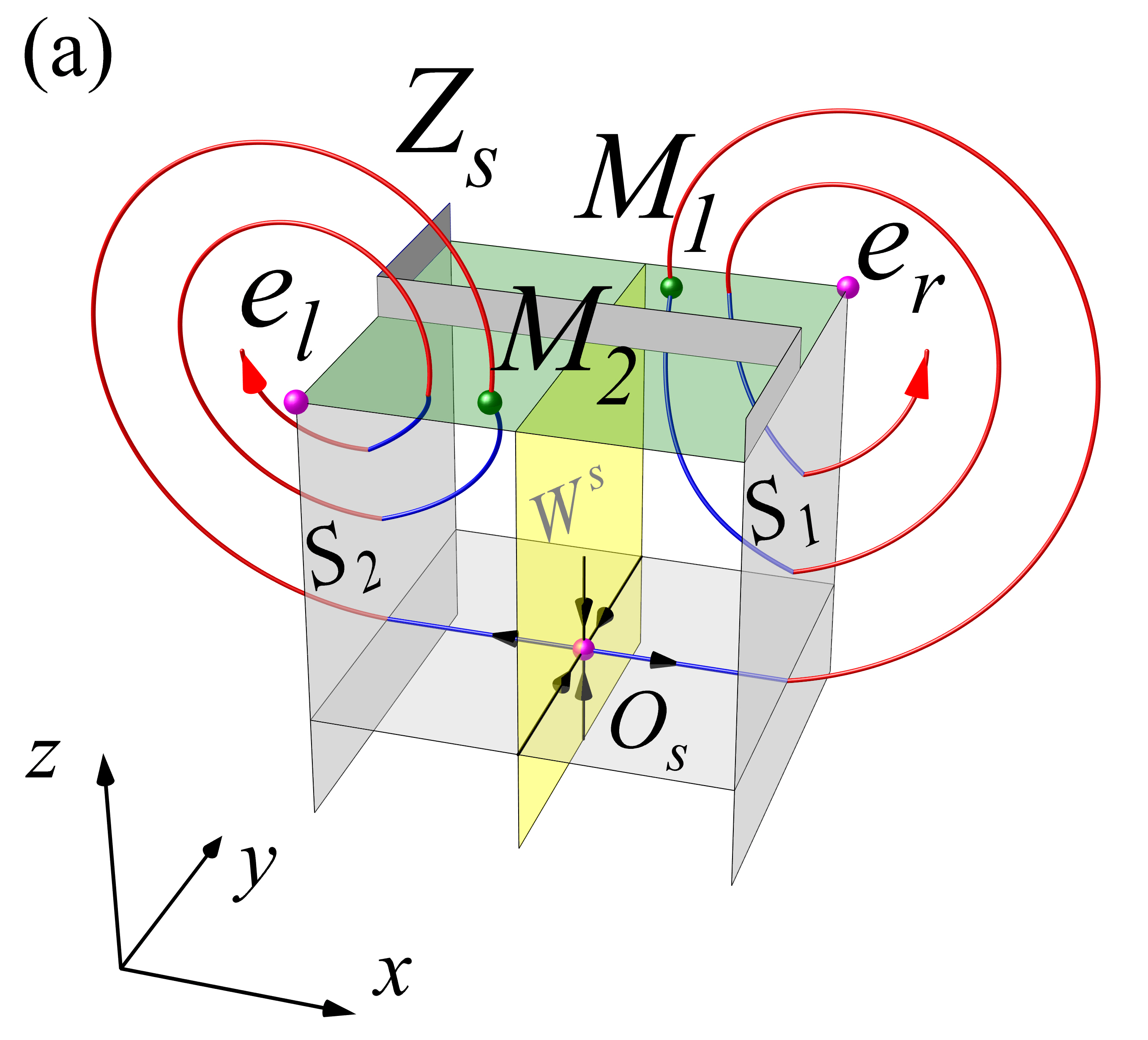}
        \includegraphics[width=0.32\textwidth]{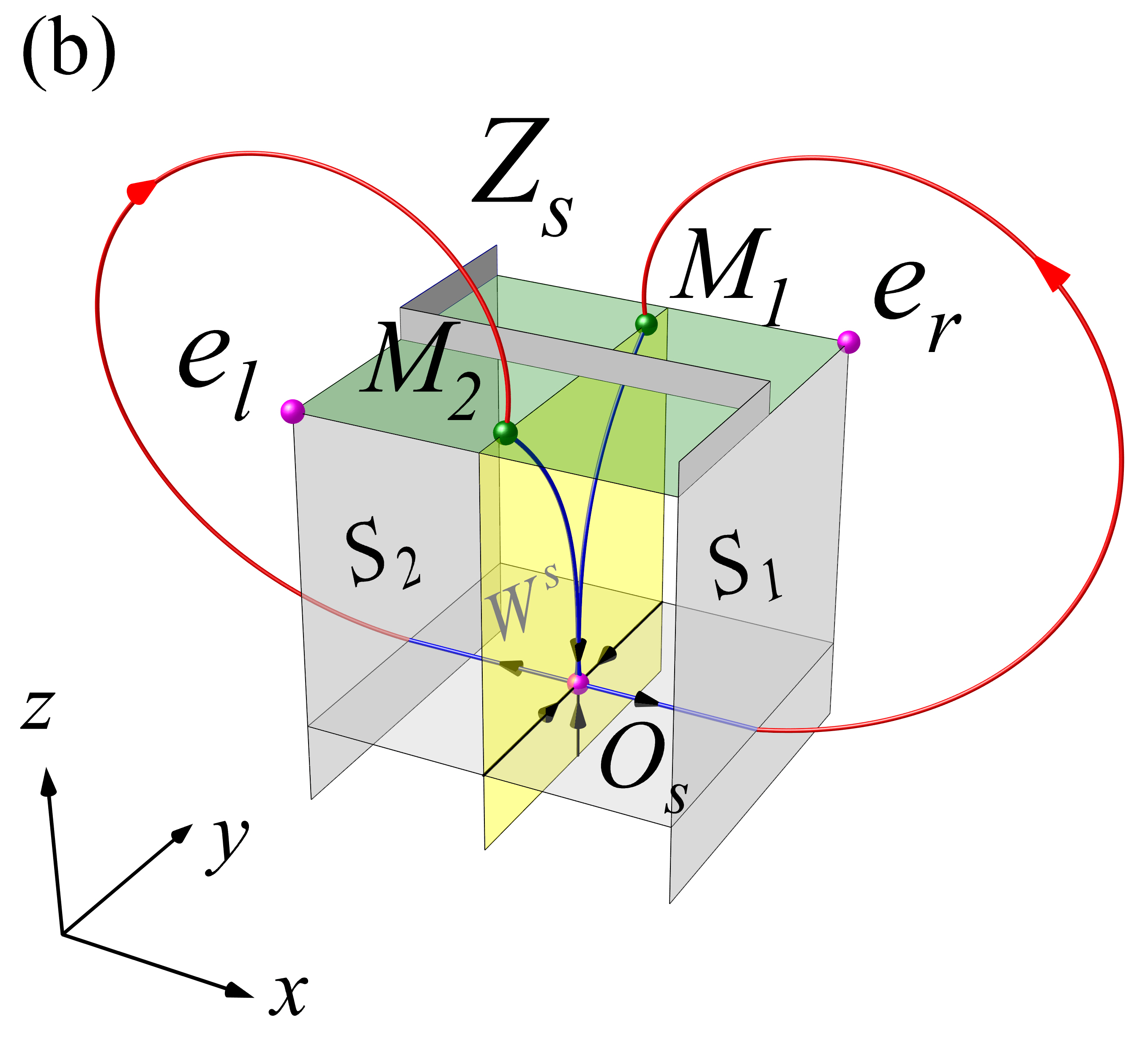}
            \includegraphics[width=0.32\textwidth]{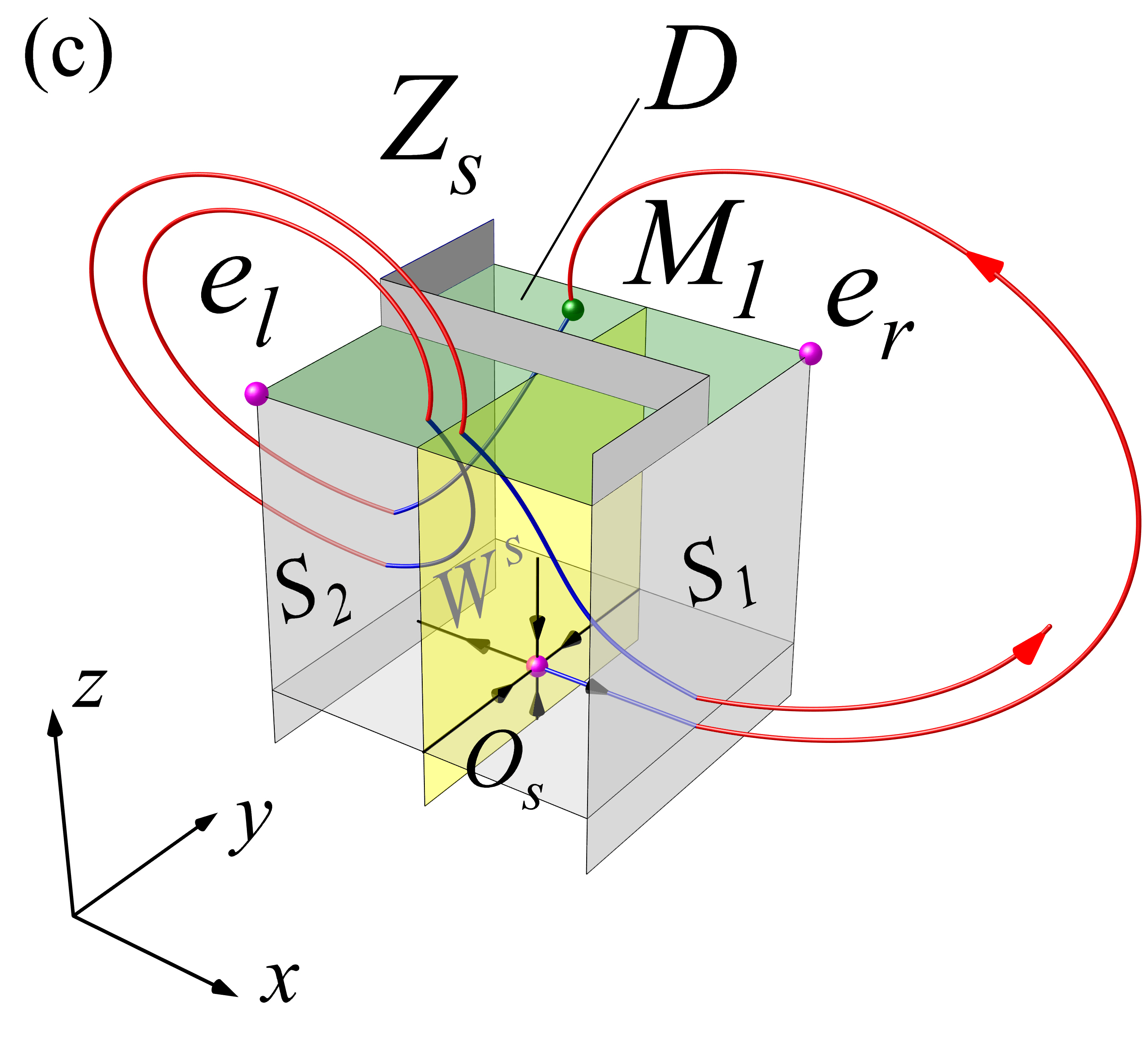}
    \caption{Formation of a homoclinic orbit in in PWL system~\eqref{sys:PWL_Lorenz} with increasing $b$ (sketch from \cite{belykh2019lorenz}).}
    \label{fig:gluing}
\end{figure}

Trajectories of PWL system~\eqref{sys:PWL_Lorenz} are composed of trajectories' pieces of linear subsystems $A_{s,l,r}$, which are “sewn” at the boundaries between regions $G_{s,l,r}$. 
In this case, a set sewn (glued together) from the trajectories of linear subsystems can be quite complex.
In particular, Fig.~\ref{fig:gluing} illustrates the process of formation of two symmetrical homoclinic orbits of the saddle $O_s$, which is a typical homoclinic set of Lorenz-type systems called a ``homoclinic butterfly''.

Similar to the Lorenz system, PWL system~\eqref{sys:PWL_Lorenz} is invariant under involution $(x,y,z)$$\rightarrow$$(-x,-y,z)$, and has three equilibrium states: 

\begin{itemize}
    \item a saddle $O_s$ at the origin with eigenvalues $s^0_1=1$, $s^0_2=-\alpha$, $s^0_3=-\nu$, and 1-D unstable manifold $W_1^u$ and 2-D stable manifold $W_2^s$;
    \item two foci $e_{l,r}( \mp 1,\mp 1, b)$ with eigenvalues $s^e_{1,2}=-\lambda\pm i\omega$, $s^e_3=-\delta$. Note that the foci, despite their stability in subsystems $A_{l,r}$, can change stability in the entire PWL system~\eqref{sys:PWL_Lorenz} due to alternating of contraction in $G_{l,r}$ and expansion in $G_s$ when passing these regions.
\end{itemize}

%\section{Poincar\'{e} return map}
In~\cite{belykh2019lorenz} it was shown that in the absence of sliding motions PWL system~\eqref{sys:PWL_Lorenz} allows to rigorously obtain 2-D piecewise-smooth Poincar\'{e} return map of the cross-section $D:\;\lbrace |x|\leq 1, |y|\leq 1, z=b\rbrace$ (green horizontal plane in Fig.~\ref{fig:gluing}) into itself.
This map $F: D\rightarrow D$ has the following explicit form
\begin{equation*}
F:\quad
    \begin{array}{ll}
     \left\{
    \begin{array}{l}
         \bar{x} = 1 - \gamma + \gamma x^{\nu}, \\
         \bar{y} = 1 - q + q x^{\alpha}y,
    \end{array}
    \right.&\textrm{for}\; x > 0,\quad\quad
    \\
    \\
    \left\{
    \begin{array}{l}
    \bar{x}=\gamma - 1 - \gamma |x|^{\nu},\\
   \bar{y} = \gamma - q - 1 + q |x|^{\alpha}y,
    \end{array}
    \right.&\textrm{for}\; x < 0,
    \end{array}
     \label{2D_map_Lorenz}
\end{equation*}
where $\gamma=be^{-\frac{3\pi\lambda}{2\omega}}$ and $q=e^{-\frac{3\pi\delta}{2\omega}}$.

The map $F$ is of a triangular form, i.e. the first equation for $x$ does not depend on the variable $y$ and is the factor equation for the slave driven second equation for $y$, which is contractive along $y$.
This fact allows one to reduce study of attaractors and bifurcations of 2-D map $F$ [and hence 3-D PWL system~\eqref{sys:PWL_Lorenz}] to that of 1-D piecewise-smooth factor map
\begin{equation}
    f:\quad \bar{x}=\left(1-\gamma+\gamma|x|^\nu\right)\textrm{sign}\;x.
    \label{map:master}
\end{equation}

The factor map~\eqref{map:master} has a discontinuity at the point $x = 0$, generated by the two-dimensional stable manifold $W_2^s$ of the saddle $O$. 
The map in the zero is defined as follows
\begin{equation*}
    \begin{array}{c}
    \bar{x} = f(\pm 0) =
    \left\{
    \begin{array}{lr}
         1 - \gamma, & \quad x =+0, \\
         \gamma - 1, &\quad x =-0.
    \end{array}
    \right.\\
    \end{array}
    \label{eq:discontinuity}
\end{equation*}

For $\gamma = 1$, the point $x = 0$ is mapped into itself, which corresponds to the birth of the two symmetrical homoclinic orbits of the saddle $O_s$ in PWL system~\eqref{sys:PWL_Lorenz}.

Like the smooth flows of Lorenz type, the attractors of PWL system~\eqref{sys:PWL_Lorenz} significantly depends on the sign of the saddle value
\begin{equation*}
    \varsigma=s_1^0+s_3^0=1-\nu,
    \label{eq:saddle_value}
\end{equation*}
controlled by the saddle index $\nu$. 
For $\nu<1$, PWL system~\eqref{sys:PWL_Lorenz} posses a strange (expanding) attractor, and for $\nu>1$, a quasi-strange attractor.

We will separately consider the two cases of $\nu<1$ and $\nu>1$, and compare the observed dynamics with the those of well-known smooth counterparts: the original Lorenz system (Sect.~\ref{sec:positive}) and the Lorenz-like Lyubimov-Zaks system (Sect.~\ref{sec:negative}).

\section{The expanding case $0<\nu<1$ \label{sec:positive}}
 
The complete study of bifurcations of PWL system~\eqref{sys:PWL_Lorenz} in the case of positive saddle value is given in the paper \cite{belykh2019lorenz} (see Theorem 3). 
Here we briefly recall the main result:

\begin{itemize}
    \item[a)] For $0<b<H_1$ [white region in Fig.~\ref{fig:bif_diag_positive}(a)], where
    \begin{equation*}
        H_1=\exp \frac{3\pi\lambda}{2\omega},
    \end{equation*}
    two stable foci $e_{l,r}$ form globally stable attractive set of PWL system \eqref{sys:PWL_Lorenz}, and a nontrivial set is absent [see the typical phase portrait in Fig.~\ref{fig:gluing}(a) and Fig.~\ref{fig:phase_port_positive}(a)].
    \item[b)] For $b=H_1$ [vertical line $H_1$ in Fig.~\ref{fig:bif_diag_positive}], there are two symmetrical homoclinic orbits of the saddle $O_s$ (the homoclinic butterfly) [see Fig.~\ref{fig:gluing}(b) and Fig.~\ref{fig:phase_port_positive}(b)].
    Arbitrary small increase of $b$ from $b=H_1$ leads to homoclinic bifurcation and to birth of two saddle cycles $C_{1,2}$, whose manifolds intersect transversally and generate a complicated Cantor set $\Omega$ of saddle orbits.
    \item[c)] For $H_1<b<b_{het}=\gamma_{het}b_h$ [gray region in Fig.~\ref{fig:bif_diag_positive}(a)], where $\gamma_{het}$ is the inverse function of
    \begin{equation*}
        \nu=1+\frac{\ln2-\ln \gamma}{\ln(\gamma-1)},
    \end{equation*}
    the stable foci $e_{r,l}$ co-exist with two symmetrical saddle cycles $C_{1,2}$ and Cantor set $\Omega$ [see Fig.~\ref{fig:phase_port_positive}(c)].
    \item[d)] For $b=b_{het}$ [the curve $b_{het}$ in Fig.~\ref{fig:bif_diag_positive}(a)], the unstable manifolds $W^u_1$ of $O_s$ fall onto the stable 2-D manifolds of $C_{1,2}$, forming two symmetrical heteroclinic orbits [see Fig.~\ref{fig:phase_port_positive}(d)]. Increasing $b$ from $b=b_{het}$ leads to a heteroclinic bifurcation, which results in the birth of a strange chaotic Lorenz-type attractor.
    \item[e)] For $b_{het}<b<b_{AH}$ [yellow region in Fig.~\ref{fig:bif_diag_positive}(a)], where
    \begin{equation*}
        b_{AH}=\nu^{-1}H_1,
    \end{equation*}
    Lorenz-type attractor co-exists with two stable foci $e_{r,l}$ [see Fig.~\ref{fig:phase_port_positive}(e)].
    \item[e-f)] The set $b=b_{AH}$ [the curve $b_{AH}$ in Fig.~\ref{fig:bif_diag_positive}(a)] corresponds to a subcritical Andronov-Hopf-like bifurcation when the saddle periodic orbit $C_1$ ($C_2$) shrinks to the stable focus $e_r$ ($e_l$) and disappear, rendering $e_r$ ($e_l$) a saddle-focus.
    \item[f)] For $b_{AH}<b<b_{cr}$ [green region in Fig.~\ref{fig:bif_diag_positive}(a)], where 
    \begin{equation*}
        b_{cr}= 2\sqrt{1+\frac{\lambda^2}{\omega^2}}
 \exp\Big\lbrace\frac{\lambda}{\omega}\big(\arctan\frac{\omega}{\lambda}+\pi\big)\Big\rbrace,
    \end{equation*}
the strange Lorenz-type attractor is a unique attractor of PWL system~\eqref{sys:PWL_Lorenz} [see Fig.~\ref{fig:phase_port_positive}(f)]. 
\end{itemize}

For $b=b_{cr}$, stable sliding motions emerge inside the attractor and destroy its chaoticity. This case was considered in the paper \cite{belykh2021sliding} in details.

\subsection{Comparison with the Lorenz system}

Now we show that these bifurcations and attractors qualitatively coincide with those of the original Lorenz system
\begin{equation*}
\begin{array}{l}
    \dot{x}=\sigma(y-x),\\
    \dot{y}=x(r-z)-y,\\
    \dot{z}=xy-\frac{8}{3}z,
\end{array}
\end{equation*}
in the ``classical'' parameter range $\sigma=10$, $r\in[10,28]$.
Let us mention the most important properties of the system for us:
\begin{itemlist}
    \item the Lorenz system is odd-symmetric $(x,y,z)$$\rightarrow$$(-x,-y,z)$;
    \item the system has a saddle $O_s$ at the origin with 1-D unstable invariant manifold and 2-D stable invariant manifold.
    For the considered parameter range, its saddle value $\varsigma$ is positive ($0<\nu<1$) 
\begin{equation*}
    \varsigma=-\frac{1}{2}(1+\sigma)+\sqrt{\frac{1}{4}(1+\sigma)^2+\sigma r}-\frac{8}{3}>0;
\end{equation*}
\item for the considered parameter range, the system has two symmetrical foci (saddle-foci after Andronov-Hopf bifurcation) $e_{r,l}\left(\pm\sqrt{\frac{8}{3}(r-1)},\pm\sqrt{\frac{8}{3}(r-1)},r-1\right)$. 
\end{itemlist}

\begin{figure}[h]
\centering
\includegraphics[width=0.8\textwidth]{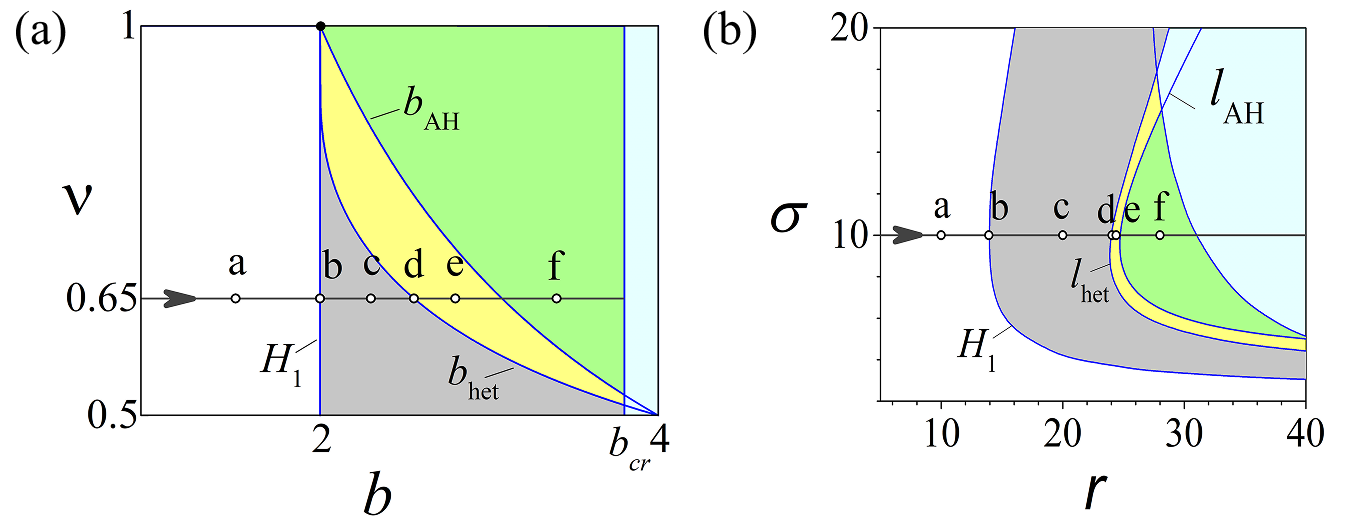}
 \caption{Comparison of bifurcation diagrams of PWL system~\eqref{sys:PWL_Lorenz} [the panel (a)] and the Lorenz system [panel (b)]. 
 The curves $b_h$ and $l_h$ are homoclinic bifurcations, $b_{het}$ and $l_{het}$ are heteroclinic bifurcations, $b_{unq}$ and $l_{AH}$ are Andronov-Hopf bifurcations. 
In both systems, the main bifurcations occur in the same order.
Lorenz strange attractor exists in yellow and green regions.  
 The route points: [the panel (a)] a$(0.65, 1.5)$, b(0.65, 2), c$(0.65, 2.3)$, d$(0.65, 2.556)$, e$(0.65, 2.8)$, f$(0.65, 3.4)$; [the panel (b)] a(10, 10), b(13.927, 10), c(20, 10), d(24.06, 10), e(24.4, 10), f(28, 10).
  The parameters $\alpha = 2$, $\nu = 0.65$, $\lambda = 0.294$, $\omega = 2$, and $\delta = 0.588$ of PWL system~\eqref{sys:PWL_Lorenz} are fixed. The line $b_{cr}=3.95548$ corresponds to destruction of the attractor via appearance of sliding motions.
 Phase portraits for corresponding route's points are in Fig.~\ref{fig:phase_port_positive}.} 
     \label{fig:bif_diag_positive}
\end{figure}

Let us compare the bifurcation diagrams of both systems, depicted in Fig.~\ref{fig:bif_diag_positive}. 
The diagram of PWL system~\eqref{sys:PWL_Lorenz} [the panel (a) in Fig.~\ref{fig:bif_diag_positive}] was proven in the paper \cite{belykh2019lorenz} and the diagram of the Lorenz system [the panel (b) in Fig.~\ref{fig:bif_diag_positive}] was obtained numerically in a large number of papers \cite{shilnikov1980appendix,sparrow1982lorenz,bykov1989boundaries,barrio2012kneadings,creaser2017finding}.

\begin{figure}[h]
\centering
\includegraphics[width=0.8\textwidth]{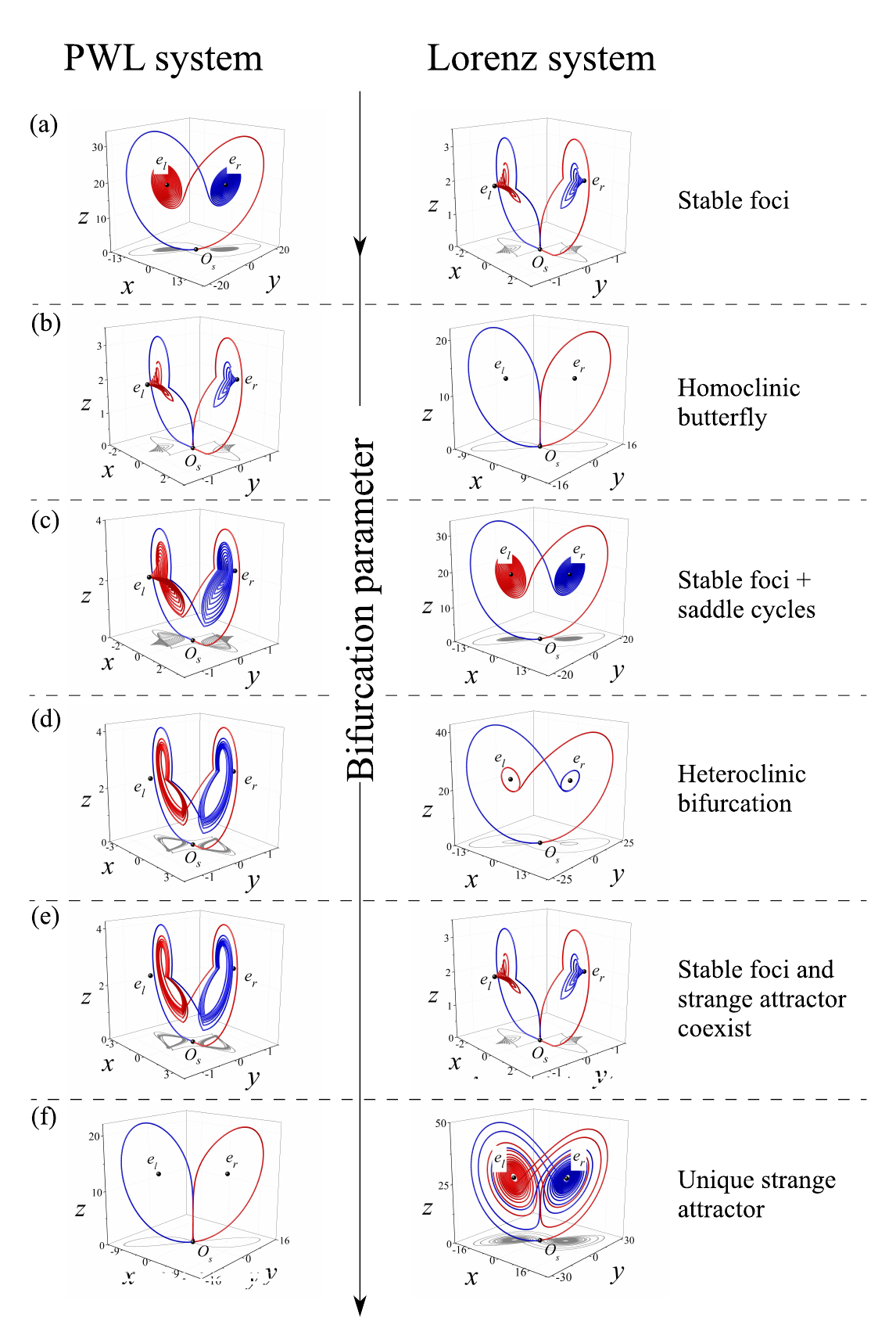}
\caption{Comparison of phase portraits of PWL system~\eqref{sys:PWL_Lorenz} and the Lorenz system when the bifurcation parameter changes along the route $\text{a}\rightarrow \text{f}$ (see diagrams in Fig.~\ref{fig:bif_diag_positive}). The rows (a)--(f) correspond to points $\text{a} - \text{f}$ in Fig.~\ref{fig:bif_diag_positive} with the same parameter values.}
\label{fig:phase_port_positive}
\end{figure}

For both diagrams we choose the route $\text{a}\rightarrow \text{f}$ of codimension 1 bifurcations (line with labeled red circles in Fig.~\ref{fig:bif_diag_positive}) leading from trivial dynamics to the birth of the strange attractor.
Moving along the route, we change the parameter responsible for the $z$-coordinate of foci $e_{l,r}$ in both systems.
They are $b$ in PWL system~\eqref{sys:PWL_Lorenz} and $r$ in original Lorenz model.
For each route point, we numerically illustrated  phase portraits of both systems (see Fig.~\ref{fig:phase_port_positive}).

We got the following comparison:
\begin{itemize}
    \item[a)]
    Both systems have two stable foci $e_{l,r}$, attracting the entire phase space except the stable manifold of the saddle $O_s$.
One-dimensional unstable manifolds $W_{1,2}^u$ of $O_s$ tend to the foci such that right branch $W_1^u\rightarrow e_r$ and $W_2^u\rightarrow e_l$ [see phase portraits in Fig.~\ref{fig:phase_port_positive}(a)].
\item[b)] 
This point belongs to the homoclinic bifurcation curve $H_1$, where $W_{1,2}^u$ fall on 2-D stable manifold $W^s$ of the saddle $O_s$ forming homoclinic butterfly [see Fig.~\ref{fig:phase_port_positive}(b)].
This bifurcation leads to the birth of a pair of saddle limit cycles $C_{1,2}$, as well as a nontrivial limiting set $\Omega$ of saddle orbits.
\item[c)] Both system have co-existing stable foci and the non-trivial limiting set, born at the bifurcation of the homoclinic butterfly.
Unstable manifolds $W_{1,2}^u$ tend to the ``opposite'' foci such that $W_{1}^u\rightarrow e_l$ and $W_{2}^u\rightarrow e_r$ [see Fig.~\ref{fig:phase_port_positive}(c)]. 
\item[d)] 
The route intersects the heteroclinic bifurcation (curves $b_{het}$ and $l_{het}$ in Fig.~\ref{fig:bif_diag_positive}) when $W_{1,2}^u$ fall on 2-D stable manifolds of the saddle orbits $C_{1,2}$. 
This bifurcation leads to the birth of the strange attractor via forming a trapping region containing orbits of saddle-type only [see Fig.~\ref{fig:phase_port_positive}(d)]. 
\item[e)] 
Stable foci $e_{l,r}$ and the strange attractor co-exist [see yellow region in Fig.~\ref{fig:bif_diag_positive} and  Fig.~\ref{fig:phase_port_positive}(e) for phase portraits].
\item[f)] 
Finally, the route intersects the curve of Andronov-Hopf bifurcation ($b_{AH}$ and $l_{AH}$ in Fig.~\ref{fig:bif_diag_positive}), where $C_{1,2}$ shrink at foci $e_{l,r}$ making them unstable, and enters the parameter region (green region in the diagrams), where the strange attractor is the unique attracting set of both systems [see Fig.~\ref{fig:phase_port_positive}(f)].  
\end{itemize}

From a comparison of the sequence of bifurcations and phase portraits of the two systems, a qualitative coincidence of dynamics follows.

\section{The contracting case $\nu>1$\label{sec:negative}}

Consider PWL system~\eqref{sys:PWL_Lorenz} in the case of negative saddle value $\varsigma<0$, which corresponds to the saddle index $\nu>1$.

For convenience, we introduce the following notions. %of a homoclinic butterfly of period $n$ and a limit cycle of type eight of period $n$.

\textbf{Denotation 1.}
We will say that a homoclinic orbit of the saddle $O$ is $n$-round if the piece of 1-D unstable manifold $W_{1,2}^u$ $n$-times intersects the Poincar\'{e} cross section $D$.

\textbf{Denotation 2.}
A pair of symmetrical $n$-round homoclinic orbits we will call \textit{$n$-round homoclinic butterfly} $B^n$.
This butterfly corresponds to a pair of $n$-periodic orbits of the factor map~\eqref{map:master} which contain the singular point $x = 0$
\begin{equation*}
    f^n(0)=0.
\end{equation*}

\textbf{Denotation 3.}
The periodic orbit of the system~\eqref{sys:PWL_Lorenz} having $n$ distinguishable intersections with the cross-section $D$ we will call the $n$-periodic limit cycle $C^n$ of figure-eight type.
This cycle corresponds to $n$-periodic orbit of the map~\eqref{map:master}, which does not contain the point $x = 0$.

\subsection{Bifurcations of the factor map for $\nu>1$ \label{subsec:bifs_1_D_nu>1}}

Attractors and bifurcations of PWL system~\eqref{sys:PWL_Lorenz} with a negative saddle value correspond to those of the 1-D factor map~\eqref{map:master} with $\nu>1$.
Let us consider its properties, taking $\gamma$ as the bifurcation parameter.

For the whole interval $\gamma\in[0,\gamma_{cr}]$ the map~\eqref{map:master} has two fixed points $e_{l,r}: x=\pm 1$, which multiplier is defined as 
\begin{equation*}
\mu_e=\gamma\nu.
\end{equation*}

For $\gamma=0$, fixed points $e_{l,r}$ are superstable due to $\mu_e=0$ and attract the interval $x\in[-1,1]$ in one iteration: $f(x)=e_l$ for $x\in[-1,0)$ and $f(x)=e_r$ for $x\in(0,1]$. 

\begin{figure}[h]  
%\vspace{-4ex} 
\centering 
 \includegraphics[width=0.3\linewidth]{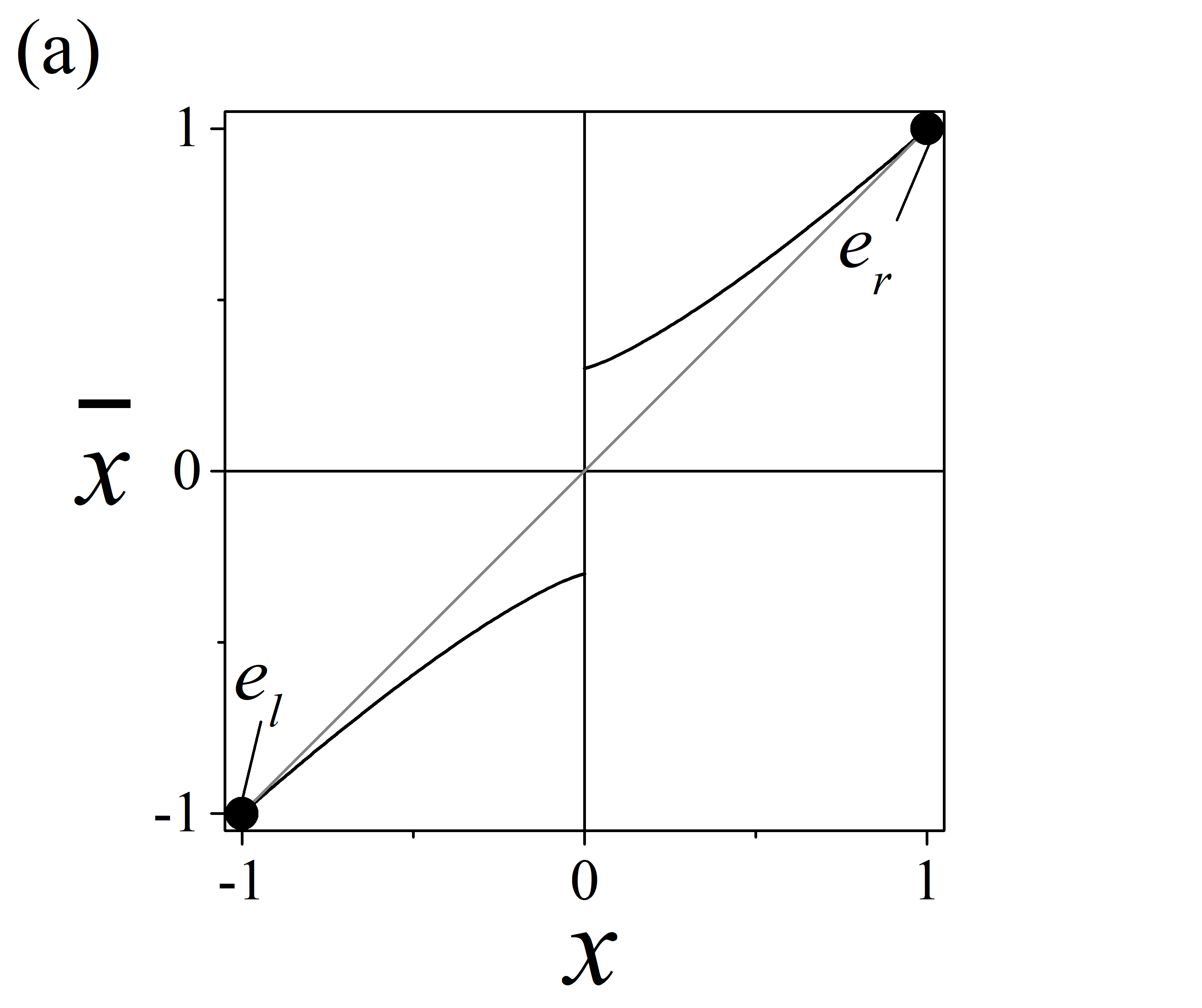} 
\includegraphics[width=0.3\linewidth]{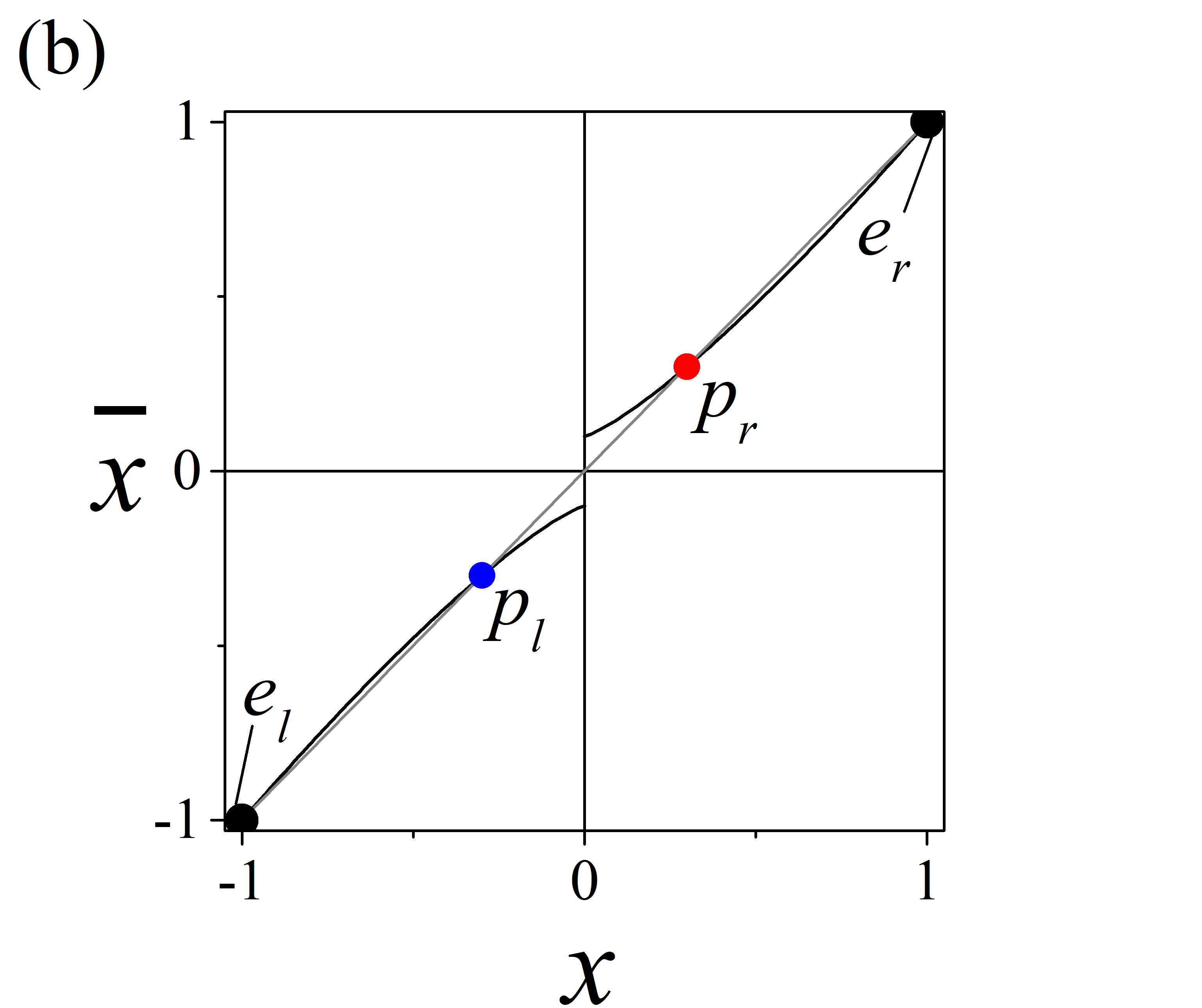}    \includegraphics[width=0.3\linewidth]{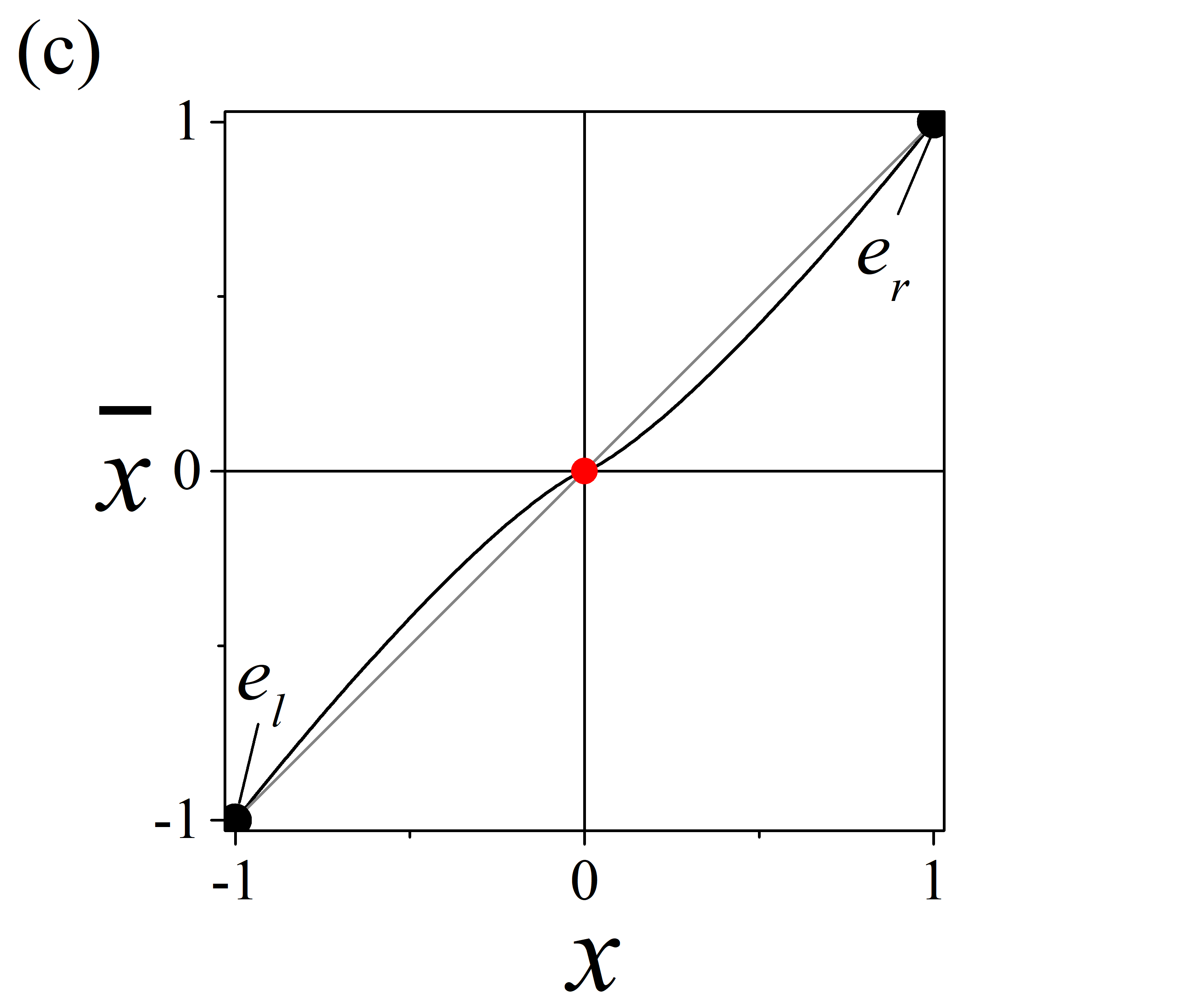} 
  \includegraphics[width=0.3\linewidth]{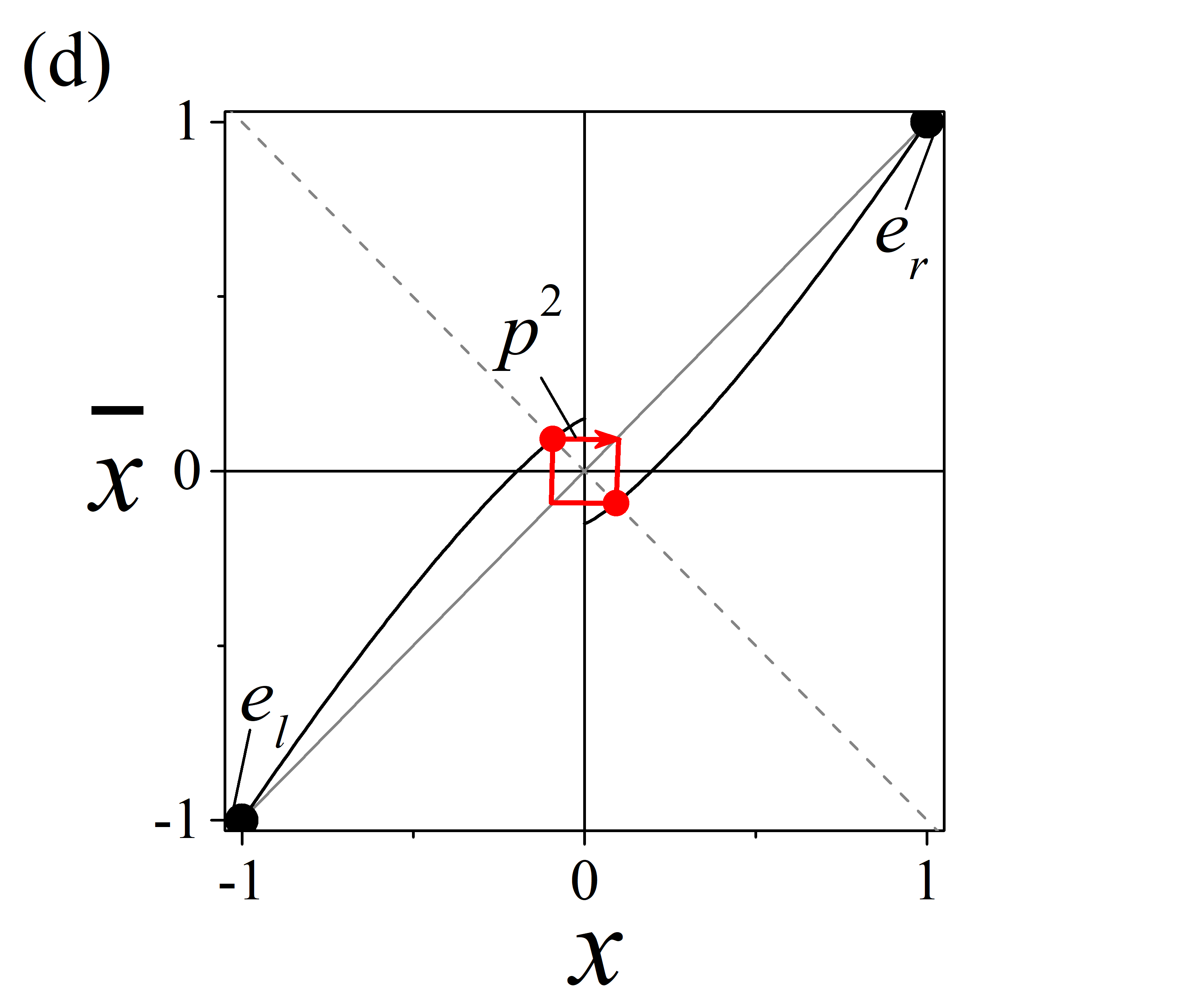}
  \includegraphics[width=0.3\linewidth]{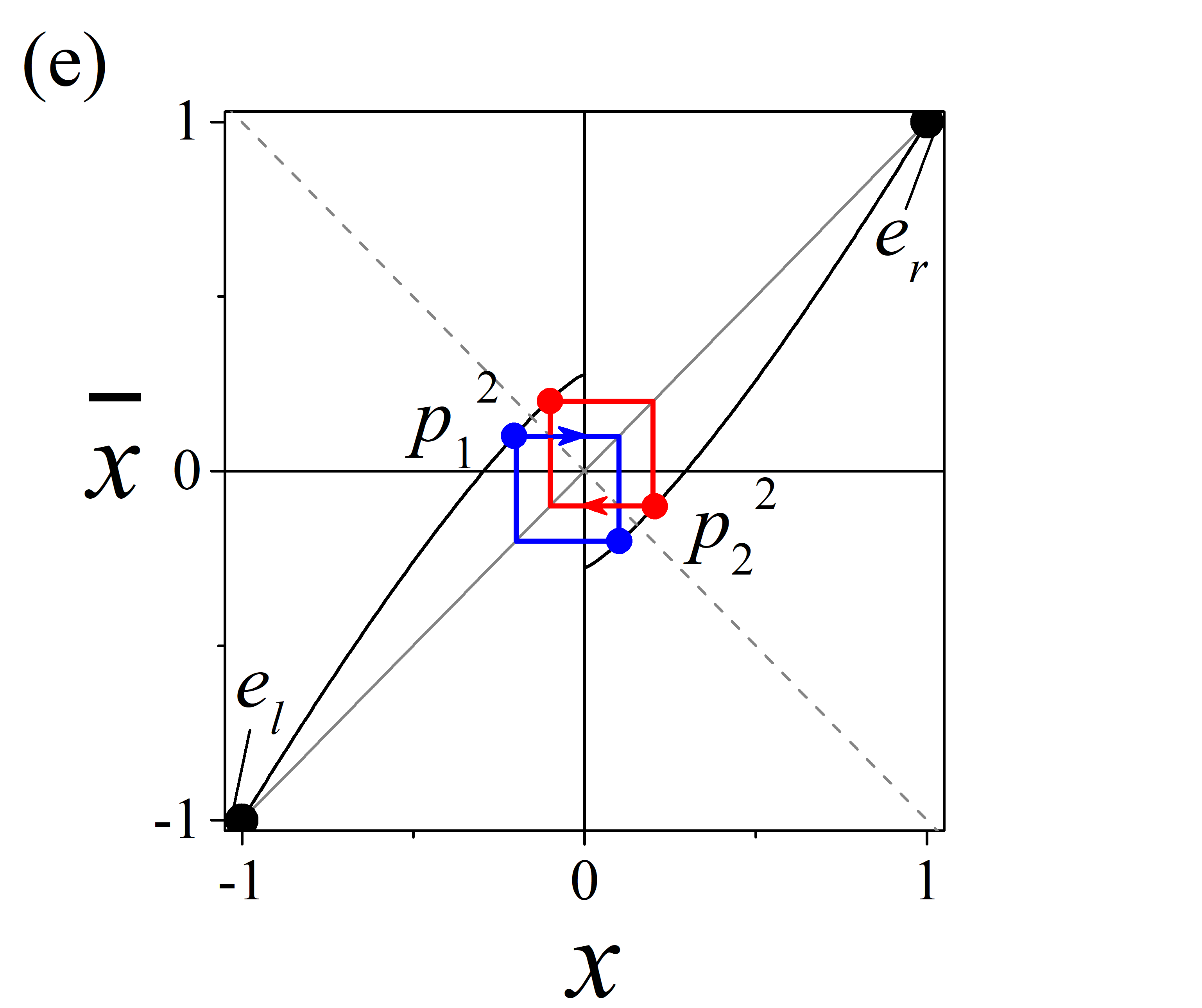}
  \includegraphics[width=0.3\linewidth]{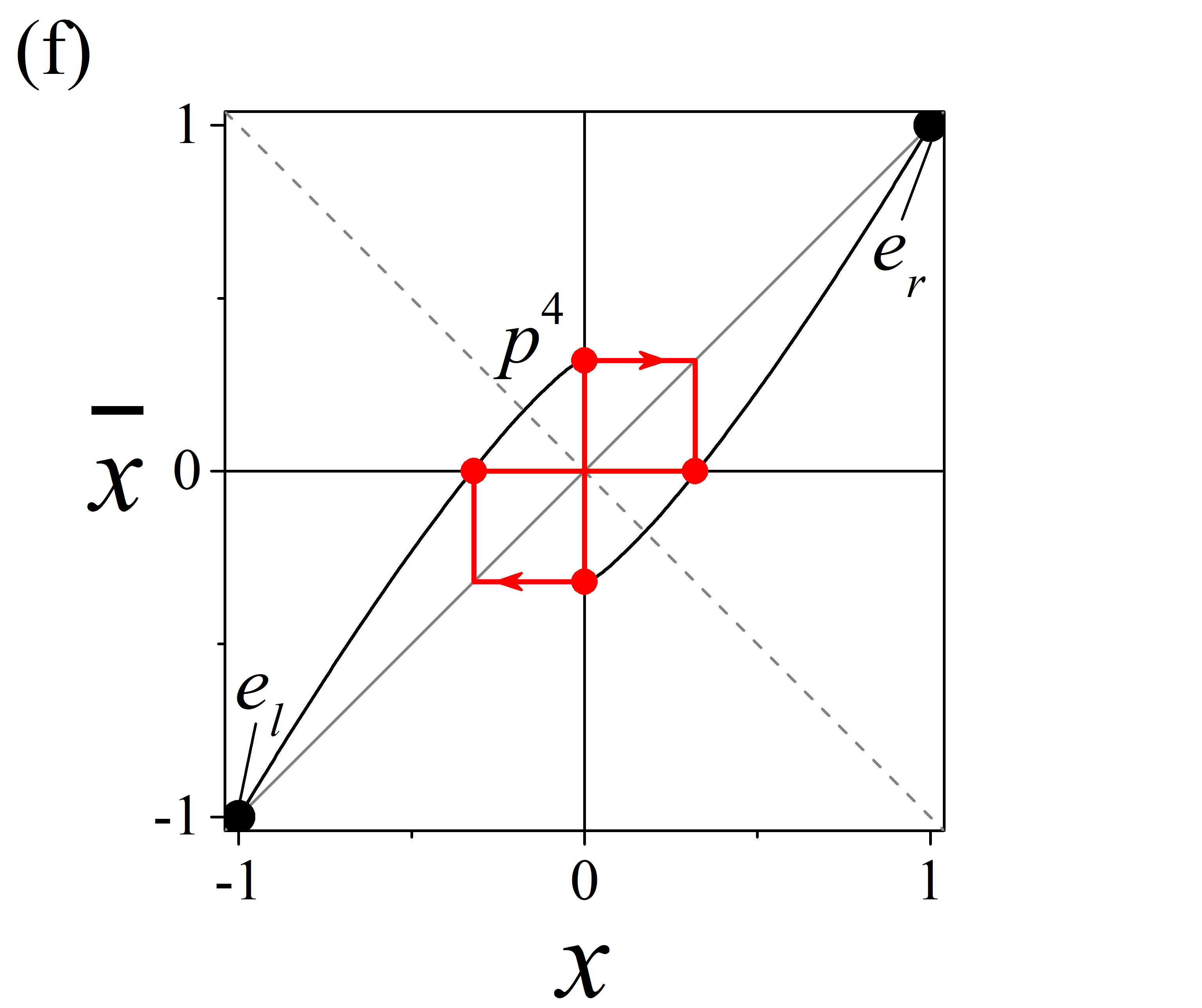}
\caption
{
    The first bifurcations of the factor map~\eqref{map:master}.
    (a) $\gamma$ = 0.7. Fixed points $e_{l,r}$ attract the entire interval.
    (b) $\gamma$ = 0.8. Attractors are fixed points $p_{l,r}$, born from $e_{l,r}$ as a result of a transcritical bifurcation.
    (c) $\gamma$ = 0.9. The points $p_{l,r}$ formed a superstable fixed point at zero.
    (d) $\gamma$ = 1. The first stable orbit $p^2$ of period 2 was born.
    (e) $\gamma$ = 1.277.Attractors are orbits $p_{1,2}^2$ of period two, born from $p^2$ as a result of a pitchfork bifurcation.
    (f) $\gamma$ = 1.325. Orbits $p_{1,2}^2$ stuck with zero, thus forming a superstable orbit $p^4$ of period 4.
}
\label{fig:map_v>1}
\end{figure}

Increasing the parameter $\gamma$ leads to the following scenario.
\begin{alphlist}[(f)]
\item[a)] 
For $0<\gamma<\gamma_{AH}$, where 
\begin{equation*}
\gamma_{AH}=\nu^{-1},
\end{equation*}
the points $e_{l,r}$ remain stable with the basins $[-1,0)$ for $e_l$ and $(0,1]$ for $e_r$ [see the graph of the map in Fig.~\ref{fig:map_v>1}(a)].

\item[b)] 
For $\gamma=\gamma_{AH}$, the multiplier $\mu_e=1$ and the points $e_{l,r}$ loose their stability via transcritical bifurcation.

Slightly increase in $\gamma$ from $\gamma=\gamma_{AH}$ leads to the emergence of a pair of stable fixed points $p_{l,r}$ from $e_{l,r}$ [red and blue dots in Fig.~\ref{fig:map_v>1}(b)].

\item[c)] 
For $\gamma=H_1$, 
where
\begin{equation*}
    H_1=1,
\end{equation*}
the points $p_{l,r}$ merge at the origin $x=0$ and form a superstable fixed point due to zero multiplier $\mu_p=0$ and has the whole interval $(-1,1)$ as a basin [see Fig.~\ref{fig:map_v>1}(c)].

\item[d)] 
A slight increase in $\gamma$ from value $\gamma=1$ leads to the birth of the stable orbit $p^2$ of period 2. 

This orbit remains stable for $1<\gamma<F_2$, where $\gamma=F_2$
is the root of the equation
\footnote{This equation was obtained from the system of two equations $1-\gamma+\gamma x^\nu=-x$, $\gamma\nu x^{\nu-1}=1$, meaning that the period 2 orbit $p^2$ has a multiplier $+1$.}
\begin{equation*}
\gamma\nu^{\nu}\left( \frac{\gamma-1}{\nu+1} \right)^{\nu-1}=1.
%\left( \frac{1}{\nu} + 1 \right)(\nu \gamma)^{\frac{1}{1-\nu}} = \gamma - 1.
    \label{eq:pitchFork2}
\end{equation*}

\item[e)] 
For $\gamma=\gamma_{F2}$, the multiplier of the orbit $p^2$ becomes equal to $+1$ and the orbit undergoes a supercritical pitchfork bifurcation, giving rise to two stable orbits $p^2_{1,2}$ [see red and blue trajectories in Fig.~\ref{fig:map_v>1}(e)]. 

\item[f)] 
For $\gamma=H_2$, where $H_2$ is the root of the equation\footnote{This equation is obtained from the condition $f^2(0)=0$. Note, that it can be rewritten in the explicit form $\nu=\log_{\gamma-1}(\frac{\gamma-1}{\gamma})$.}
\begin{equation}
    \gamma(\gamma-1)^{\nu-1}=1,
    \label{eq:homocl_2}
\end{equation}
two period 2 orbits $p^2_1$ and $p^2_2$ merge at the $x=0$ [see Fig.~\ref{fig:map_v>1}(f)].

With a slight increase in the parameter $\gamma$ from $\gamma=\gamma_{H_2}$, a stable orbit $p^4$ of period 4 is born.
%, which corresponds to the formation of the attracting homoclinic butterfly $B^2$ in the system~\eqref{sys:PWL_Lorenz}.
\end{alphlist}

A further increase in $\gamma$ leads to the second pitchfork bifurcation $F_4$, when the orbit $p_4$ loses stability and two $p^4_{1,2}$  are born from it. 

Thus, we have a cascade of alternating homoclinic butterfly bifurcations $H_{2^n}$, defined as
\begin{equation*}
    f^{2^n}(0)=0,
\end{equation*}
and pitchfork bifurcations $F_{2^n}$, defined as
\begin{equation*}
    f^{2^n}(x)-x=0, \quad \frac{d}{dx} f^{2^n}(x)=1,
\end{equation*}
for $n\rightarrow\infty$, which leads to the birth of a quasi-strange attractor beyond the limiting curve $H_\infty$:
\begin{equation}
\encircle{$p_{l,r}$} \xrightarrow{H_1} 
\encircle{$p^{2}$} \xrightarrow{F_2} 
\encircle{$p_{1,2}^2$} \xrightarrow{H_2} 
\encircle{$p^{4}$} \xrightarrow{F_4}\cdots
 \xrightarrow{H_{2^n}} \encircle{$p^{2^{n+1}}$}
\xrightarrow{F_{2^{n+1}}}\encircle{$p^{2^{n+1}}_{1,2}$}
\xrightarrow{H_{2^{n+1}}}
\cdots 
\xrightarrow{H_{\infty}}
\encircle{$p^{\infty}$} 
\rightarrow \cdots .
\label{eq:diagram_1d} 
\end{equation}

After infinite period orbit $p^{\infty}$ is born, there is a cascade of appearing of the orbits of all possible even periods, which ends with the emerging of the first odd period $2n+1$ orbit as a result of a homoclinic bifurcation $H_{p2}$ of period 2 orbit.

A complete study of this bifurcation set of the map~\eqref{fig:map_v>1} is out of scope of this paper. 
However, according to the results on 1-D unimodal and Lorenz maps \cite{Mira1987,belykh1997comments,Malkin2003,anuvsic2020topological,belykh20241d}, one can state the existence of self embeddings similar to Mira's ``box-within-a-box'' structure, Sharkovsky ordering, a cascade of period adding bifurcations and other basic features of such maps.

% \begin{equation}    C_{1,2}\xrightarrow{H_1}C^2\xrightarrow{F_2}C^2_{1,2}\xrightarrow{H_2}C^4\xrightarrow{F_4}\ldots
% \end{equation}

% \begin{equation}
% \begin{array}{rc}
%    \text{Homoclinic/pitchfork doubling cascade:}  &  \encircle{$p_{l,r}$} \xrightarrow{H_1} 
% \encircle{$p_{2}$} \xrightarrow{F_2} 
% \encircle{$p_{1,2}^2$} \xrightarrow{H_2} 
% \encircle{$p_{4}$} \xrightarrow{F_4}\cdots
% \xrightarrow{H_{2^\infty}}
% \encircle{$\mathcal{A}_{F}$} 
% \rightarrow \\
%    \text{homoclinic bif-s of even periodic orbit:}  & \cdots\xrightarrow{H_{p_{2}}}
% \encircle{$p_{2n+1}$}\rightarrow\cdots
% \end{array}
% \label{eq:diagram_1d}
% \end{equation}

% \begin{equation}
% \begin{tikzcd}[row sep=tiny, column sep=small]
% \encircle{$p_l$} \arrow[dr, "H_1"] && \encircle{$p_2^1 $} \arrow[dr, "H_2"] && \encircle{$p_4^1$} \arrow[dr, "H_4"] &&\encircle{$p_{2^n}^1$} \arrow[dr, "H_4"]&&& \\
% & \encircle{$p_2$} \arrow[ur, "F_2"] \arrow[dr]  && \encircle{$p_4$}  \arrow[ur, "F_4"]\arrow[dr] && \cdots \arrow[ur, "F_{2^n}"]\arrow[dr] && \encircle{$\mathcal{A}_{F}$} \arrow[r] & \cdots \arrow[r, "H_{p2}"] & \encircle{$p_{2n+1}$} \arrow[r]&\cdots \\
% \encircle{$p_r$} \arrow[ur] && \encircle{$p_2^2$} \arrow[ur] && \encircle{$p_4^2$} \arrow[ur] &&\encircle{$p_{2^n}^2$} \arrow[ur]&&&
% \end{tikzcd}
% \label{eq:diagram_1d}
% \end{equation}

Starting from the value $\gamma=\gamma_{cr}<2$ (from corresponding value $b=b_{cr}$) sliding motions appear in PWL system~\eqref{sys:PWL_Lorenz} and 1-D factor map~\eqref{map:master} no longer corresponds to the flow and requires changes to its form to account sliding motions \cite{belykh2021sliding}. 
In this paper we does not consider this significantly non-smooth case. 

\subsection{Comparison with Lorenz-Lyubimov-Zaks system}

In this subsection we show that PWL system~\eqref{sys:PWL_Lorenz} with a negative saddle value reconstructs well the bifurcation properties of nonlinear smooth Lorenz-Lyubimov-Zaks (LLZ) system \cite{lyubimov1983two}.

LLZ system is an extension of the Lorenz system and has the form
\begin{equation}
 \begin{array}{l}
 \dot{x} = -\sigma x + \sigma y + \sigma y D (z - r), \\
 \dot{y}  = r x - y - x z,\\
 \dot{z}  = x y - \beta z,
 \end{array}
  \label{sys:Lyubimov}
\end{equation}
where $\sigma$, $\beta$, $r$, $D$ are positive parameters.
Note, that for $D=0$, this system becomes original Lorenz model.

A detailed numerical study of LLZ system was presented in \cite{kazakov2021bifurcations}.
Using these results (thanks to the agreement of Alexei Kazakov), we will briefly present the properties of the system and compare the bifurcation routes in PWS and LLZ systems.

\subsubsection{Basic properties of LLZ system \label{subsec:LLZ_properties}}

The system~\eqref{sys:Lyubimov} inherits from the Lorenz system odd symmetry 
\begin{equation*}
    (x, y, z) \rightarrow (-x, -y, z).
    \label{symmetry_Lyub}
\end{equation*}

Depending on the parameters, the system has from one to five equilibrium states.
In this paper, we consider Lorenz-like partitions of the phase space and the corresponding parameter regions, where the system has three equilibria
\begin{equation*}
O_s(0, 0, 0),\quad
e_{l}(-Q,QR, r+R),\quad
e_{r}(Q,-QR, r+R),
\end{equation*}
where $Q=\sqrt{\frac{\beta}{2}}\sqrt{r(\sqrt{1-4D}+1)-2}$ and $R=\frac{\sqrt{1-4D}-1}{2D}$.

As in the Lorenz system, equillibria $e_{l,r}$ are born stable from $O_s$ as a result of a supercritical pitchfork bifurcation.\footnote{For $l_{pf}>0$ and $r>2$ it subcritical branch gives rise two more saddles and the system will have 5 equilibrium states [see the curve $l_{pf}$ in Fig.~\ref{fig:2D_diags_negative}(b)].} 
\begin{equation*}
    l_{pf}(D,r)\triangleq Dr^2-r+1=0,\quad r<2.
\end{equation*}

For $l_{pf}(D,r)<0$ and $r<0$, the point $O_s$ becomes a saddle with eigenvalues
\begin{equation*}
    s_{1,2}^0=-\frac{1}{2}(1+\sigma)\pm\frac{1}{2}\sqrt{(1+\sigma)^2-4\sigma(Dr^2-r+1)},\quad s_3^0=-b,
    \label{eq:LZ_eigenvalues}
\end{equation*}
thus, similarly to the case of Lorenz system, $O_s$ has 1-D unstable manifold $W^u_1$ and 2-D stable manifold $W^u_2$.
Under the condition $s_2^0<-\beta$, 
an eigenvector $V_{s_3^0}=(0,0,1)$, corresponding to the eigenvalue $s_3^0$, defines leading direction.
Thus the saddle value $\varsigma$ of $O_s$ is defined as follows
\begin{equation*}
\varsigma\triangleq s_1^0+s_3^0=-\frac{1}{2}(1+\sigma)+\frac{1}{2}\sqrt{(1+\sigma)^2-4\sigma(Dr^2-r+1)}-\beta.
\end{equation*}

The condition of $\varsigma=0$, providing $\nu=1$, reads  
\begin{equation*}
    l_{\nu=1}(D,r)\triangleq\beta^2+(1+\sigma)\beta+\sigma(Dr^2-r+1)=0,
\end{equation*}
and $\varsigma$ is negative for $l_{\nu=1}(D,r)>0$ [see the region to the right of the curve $ l_{\nu=1}$ in Fig.~\ref{fig:2D_diags_negative}(b)].
  
Characteristic equation for $e_{l,r}$ reads
\begin{equation*}
s^3+(\sigma+\beta+1)s^2+\beta(\sigma+r+a_1)s+2\sigma \beta(r-1+a_2)=0.
%s^3+(\sigma+b+1)s^2+b[\sigma+r+DR(r+\sigma(r+R))]s+2\sigma b(r-1+rDR)(R+2)=0.
\end{equation*}
where $a_1=(r+\sigma R+\sigma r)RD$, $a_2=(R+2r)RD$.
Due to coordinates of equillibria $e_{l,r}$ and using Rauth-Hurwitz criterion we obtain, that $e_{l,r}$ exist and are stable in the following range of parameters  
\begin{equation*}
    2\left({\sqrt{1-4D}+1}\right)^{-1}<r<l_{AH}(D),    
\label{eq:foci_stab_condition}
\end{equation*}
where 
\begin{equation}
r=l_{AH}(D)\triangleq\frac{\sigma\left[\sigma+\beta+3+(\sigma+\beta-1)R^2D\right]}{\sigma-\beta-1-\left[(\sigma-1)^2+\beta(\sigma+1)\right]RD}
\label{eq:AH_LZ}
\end{equation}
is the Andronov-Hopf bifurcation of foci $e_{l,r}$.
For $D=0$, the condition $r=l_{AH}(0)$ becomes the well-known subcritical Andronov-Hopf bifurcation in the Lorenz system. 

In this paper, we have chosen a route in the parameters plane $(D,r)$ on which Andronov-Hopf bifurcation is supercritical, i.e. when crossing the curve $l_{AH}$, the foci $e_{l,r}$ lose stability and stable limit cycles are born from them [see curve $l_{AH}$ in Fig.~\ref{fig:2D_diags_negative}(b)].

\subsubsection{The comparison}

Let us trace the correspondence between the bifurcations and attractors of PWL system~\eqref{sys:PWL_Lorenz} and LLZ system~\eqref{sys:Lyubimov}. 

We can easily obtain the bifurcation set of PWL system~\eqref{sys:PWL_Lorenz} from those of the factor map~\eqref{map:master} by simply recalculating the parameter\footnote{For the correspondence of orbit types see Denotations 1--3}
\begin{equation*}
    b=\gamma e^{\frac{3\pi\lambda}{2\omega}}.
\end{equation*}

\begin{figure}[!ht]
  \centering
    \includegraphics[width=0.49\textwidth]{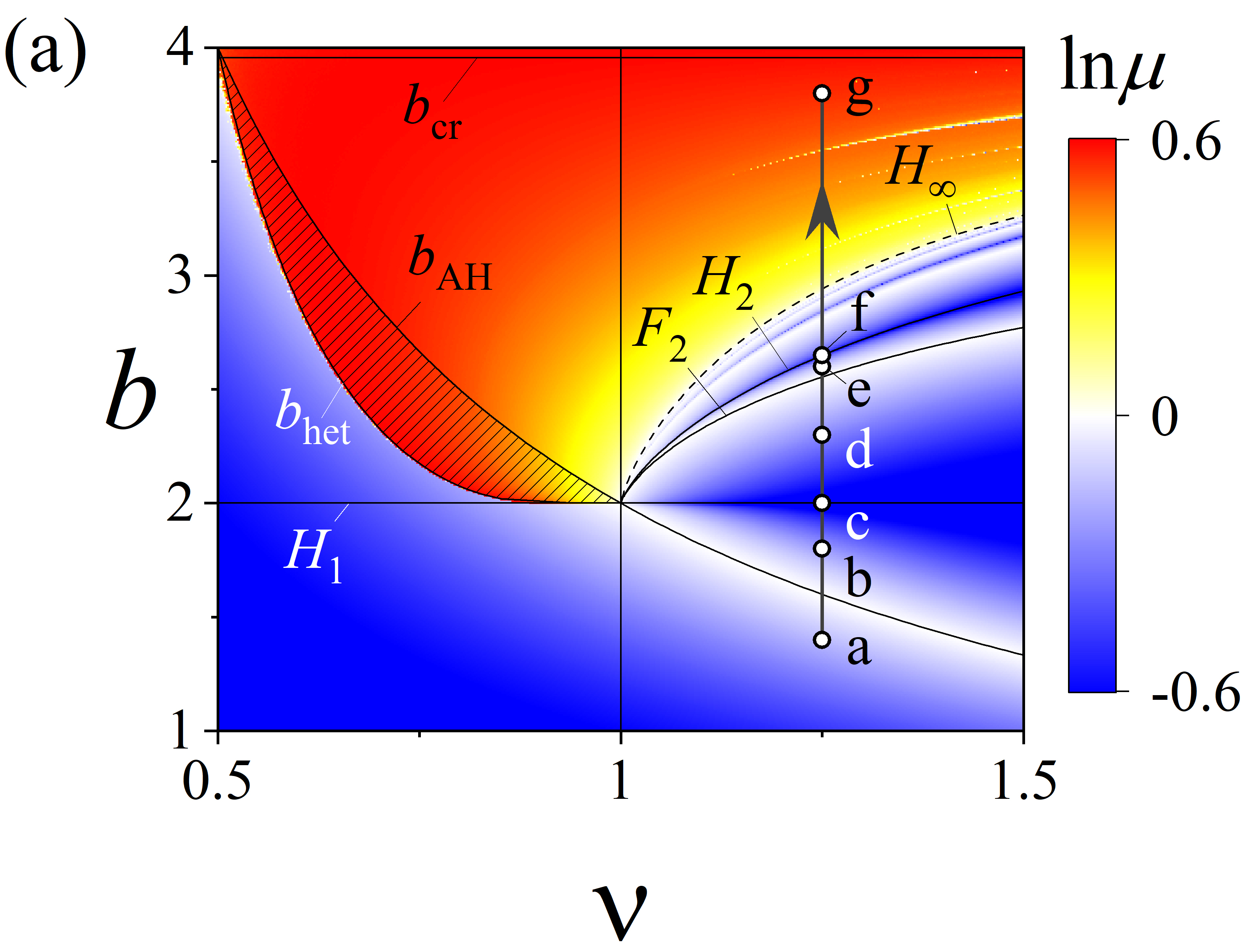}
    \includegraphics[width=0.49\textwidth]{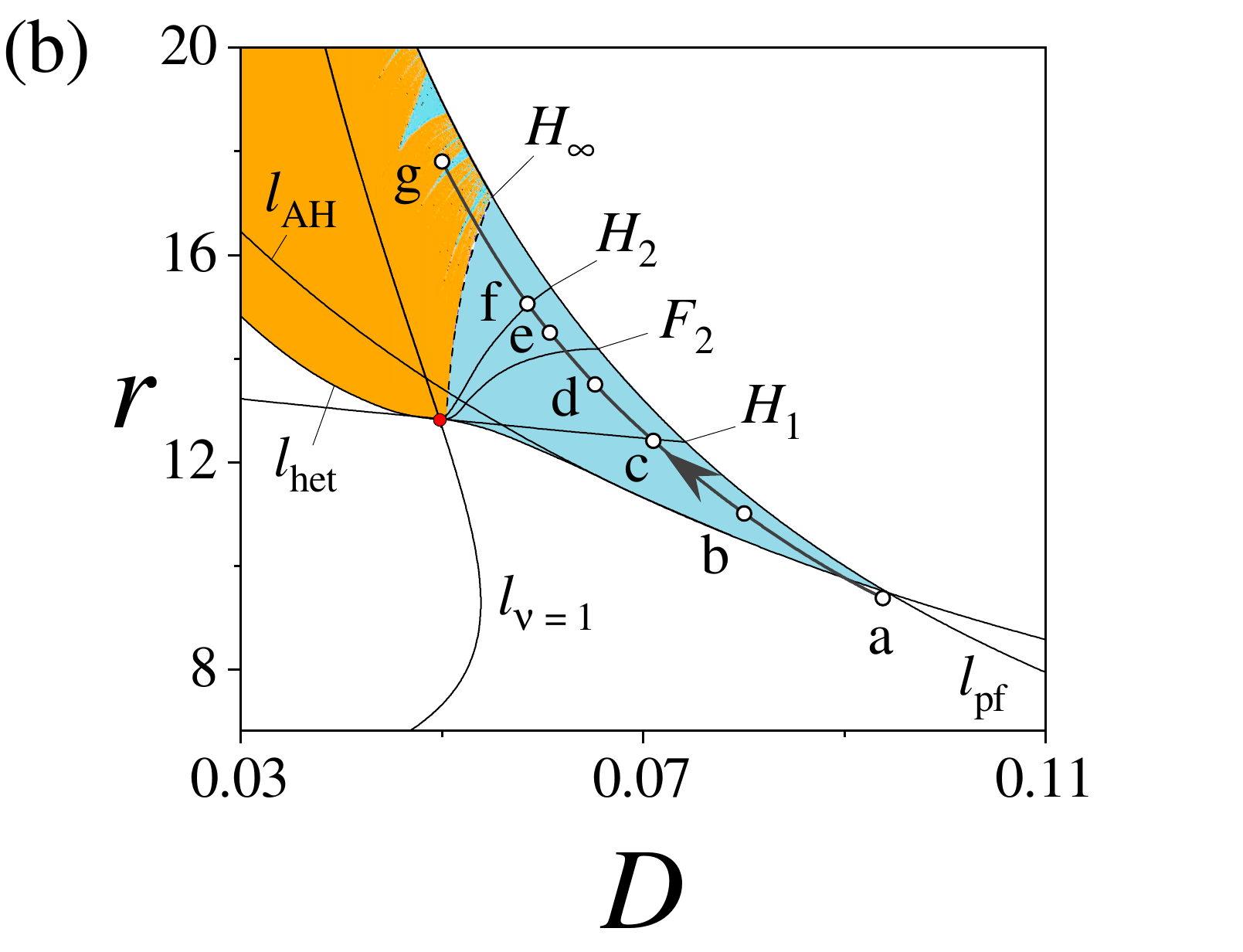}
    \caption{Comparison of bifurcation routes $\text{a}-\text{f}$ for $\nu>1$. 
    (a) Bifurcation diagram of PWL system~\eqref{sys:PWL_Lorenz}. The saddle index is fixed $\nu=1.25$. The color shows the logarithm of the factor map~\eqref{map:master} attractor multiplier $\ln\mu$. The dashed area corresponds to coexistence of stable fixed points $e_{l,r}$ and the strange attractor. The points: 
    $\text{a}(1.25, 1.4)$, 
    $\text{b}(1.25,1.8)$, 
    $\text{c}(1.25,2)$, 
    $\text{d}(1.25,2.3)$, 
    $\text{e}(1.25,2.6)$, 
    $\text{f}(1.25,2.65)$,
    $\text{g}(1.25,3.8)$. Other parameters are the same as in Fig.~\eqref{fig:bif_diag_positive}.
    (b) LLZ system~\eqref{sys:Lyubimov} (a sketch from~\cite{kazakov2021bifurcations}). Orange stands for maximal positive Lyapunov exponent and blue for zero. The route is defined by $r\approx\frac{0.881}{D}$. The points: $\text{a}(0.0938, 9.38)$, 
    $\text{b}(0.0915,9.64)$, 
    $\text{c}(0.08,11.01)$, 
    $\text{d}(0.0711,12.41)$,
    $\text{e}(0.0624,14.2)$, 
    $\text{f}(0.0586,15.1)$ ,
    $\text{g}(0.05,17.8)$. The parameters: $\sigma=10$, $\beta=\frac{8}{3}$. 
    For both systems The routes intersect the same bifurcations: $b_{AH}$ and $l_{AH}$ - Andronov-Hopf bifurcation, $H_1$ - the first (one-round) homoclinic butterfly, $F_2$ - pitchfork bifurcation, $H_2$ - two-round homoclinic butterfly, $H_\infty$ - limiting curve of quasi-strange Lorenz-type attractor birth.
    In both systems, there are periodic windows in the region above the curve $H_{\infty}$.
    The routes points are the same as in Fig.~\ref{fig:map_v>1}. 
    The point $\text{g}$ is the quasi-strange attractor.
    The corresponding phase portraits see in Fig.~\ref{fig:PWL_Lubimov_compare} and Fig.~\ref{fig:Atractor_Lorenz_like} for $\text{g}$.}
    \label{fig:2D_diags_negative}
\end{figure}

In addition to the analytically obtained bifurcation curves of the factor map~\eqref{map:master}, we numerically obtain a $(\nu,b)$-diagram of multiplier $\mu$ of its attractors [see Fig.~\ref{fig:2D_diags_negative}(a)].

The partition of the parameter space of LLZ system~\eqref{sys:Lyubimov} is given in subsection~\ref{subsec:LLZ_properties} and taken from the paper~\cite{kazakov2021bifurcations}.

For both systems we choose bifurcation routes passing through codimension 1 bifurcations in the regions with a negative saddle value $\nu>1$ (see Fig.~\ref{fig:2D_diags_negative}).
These routes are defined as: (i) $\nu=1.25$, $b\in[1.4, 3.8]$ for PWL system~\eqref{sys:PWL_Lorenz}; 
(ii) $r=\frac{0.881}{D}$, below the subcritical branch of the pitchfork bifurcation $l_{pf}$ for LLZ system~\eqref{sys:Lyubimov}.

% The fixed points $e_{l,r}$ of the map~\eqref{map:master} with their stability properties correspond to the foci $e_{l,r}(\pm 1, \pm 1, b)$ of the system~\eqref{sys:PWL_Lorenz}. 

Moving along these routes, we obtained the following bifurcation sequence.
\begin{itemize}
    \item[a)] As in the case of $\nu<1$, 
        both systems have two stable foci $e_{l,r}$, attracting the entire phase space except the stable manifold of the saddle $O_s$ [see phase portraits in Fig.~\ref{fig:PWL_Lubimov_compare}(a)].
     \item[b)] The system have two symmetric stable limit cycles $C_{1,2}$, born from foci $e_{l,r}$ via supercritical Andronov-Hopf bifurcation [curve $b_{AH}$ for PWL and $l_{AH}$ for LLZ in Fig.~\ref{fig:2D_diags_negative}].
$C_{1,2}$ correspond to fixed points $p_{l,r}$ of the factor map~\eqref{map:master} [see portraits in Fig.~\ref{fig:PWL_Lubimov_compare}(b)].
    \item[c)] The systems undergo the first homoclinic bifurcation $H_1$, when the cycles $C_{1,2}$ collide with the stable manifold $W^{s}_2$ of the saddle $O_s$ and form one-round homoclinic butterfly $B^1$, which attracts all the phase space [see Fig.~\ref{fig:PWL_Lubimov_compare}(c)]. 
    \item[d)] Both systems have figure-eight 2-periodic stable limit cycle $C^2$ born from $B^1$. The cycle $C^2$ remains stable until its pitchfork bifurcation $F_2$ [see Fig.~\ref{fig:PWL_Lubimov_compare}(d)].
    \item[e)] At the bifurcation $F_2$, cycle $C^2$ loses stability and two stable symmetric figure-eight limit cycles $C^2_{1,2}$ are born from it [see Fig.~\ref{fig:PWL_Lubimov_compare}(e)].
    \item[f)] The systems undergo the second homoclinic bifurcation $H_2$, at which the cycles $C^2_{1,2}$ stick together at the saddle $O_s$ and form 2-periodic homoclinic butterfly $B^2$ [see Fig.~\ref{fig:PWL_Lubimov_compare}(f)].
    Similar to case of the point c, $H^2$ leads to birth of 4-periodic limit cycle $C^4$, stable until pitchfork bifurcation $F_4$.
\end{itemize}

Then the pair of homoclinic $H_{2^n}$ and pitchfork $F_{2^{n}}$ bifurcations forms an infinite cascade, the limit of which is the appearance of infinite set of even $\infty$-periodic orbits $C^{\infty}$.
Beyond this limit quasi-strange Lorenz-type attractor is born.

Thus, for $\nu>1$, PWL system~\eqref{sys:PWL_Lorenz} and LLZ system~\eqref{sys:Lyubimov} undergo the same bifurcation sequence
\begin{equation}
\encircle{$C_{1,2}$} \xrightarrow{H_1} 
\encircle{$C^{2}$} \xrightarrow{F_2} 
\encircle{$C_{1,2}^2$} \xrightarrow{H_2} 
\encircle{$C^{4}$} \xrightarrow{F_4}\cdots
 \xrightarrow{H_{2^n}} \encircle{$C^{2^{n+1}}$}
\xrightarrow{F_{2^{n+1}}}\encircle{$C^{2^{n+1}}_{1,2}$}
\xrightarrow{H_{2^{n+1}}}
\cdots  
\xrightarrow{H_{\infty}}
\encircle{$C^{\infty}$} 
\rightarrow \cdots .
\label{eq:diagram_3d} 
\end{equation}

As was noted for the 1-D factor map~\eqref{map:master}, the birth of the cycle $C^{\infty}$ in a piecewise-linear system~\eqref{sys:PWL_Lorenz} is followed by a complex bifurcation structure, similar to the Mira ``box-within-a-box'' structure.
This structure contains a bifurcation of the appearance of limit cycles $C^{2n+1}$ of an odd period, the Sharkovsky ordering, a cascade of period adding and other bifurcation phenomena similar to those of unimodal maps. 

We trace the similarity of such a structure in both systems numerically using 1-D extended bifurcations diagrams (see Fig~\ref{fig:bif_diag_positive}).  

For PWL system~\eqref{sys:PWL_Lorenz}, the Poincar\'{e} cross-section is the rectangle $D$ (green plane in Fig.~\ref{fig:gluing})
\begin{equation*}
    D:\quad |x|\leq 1, |y|\leq 1, z=b.
\end{equation*}

\begin{figure}[ht]
  \centering   \includegraphics[width=0.47\textwidth]{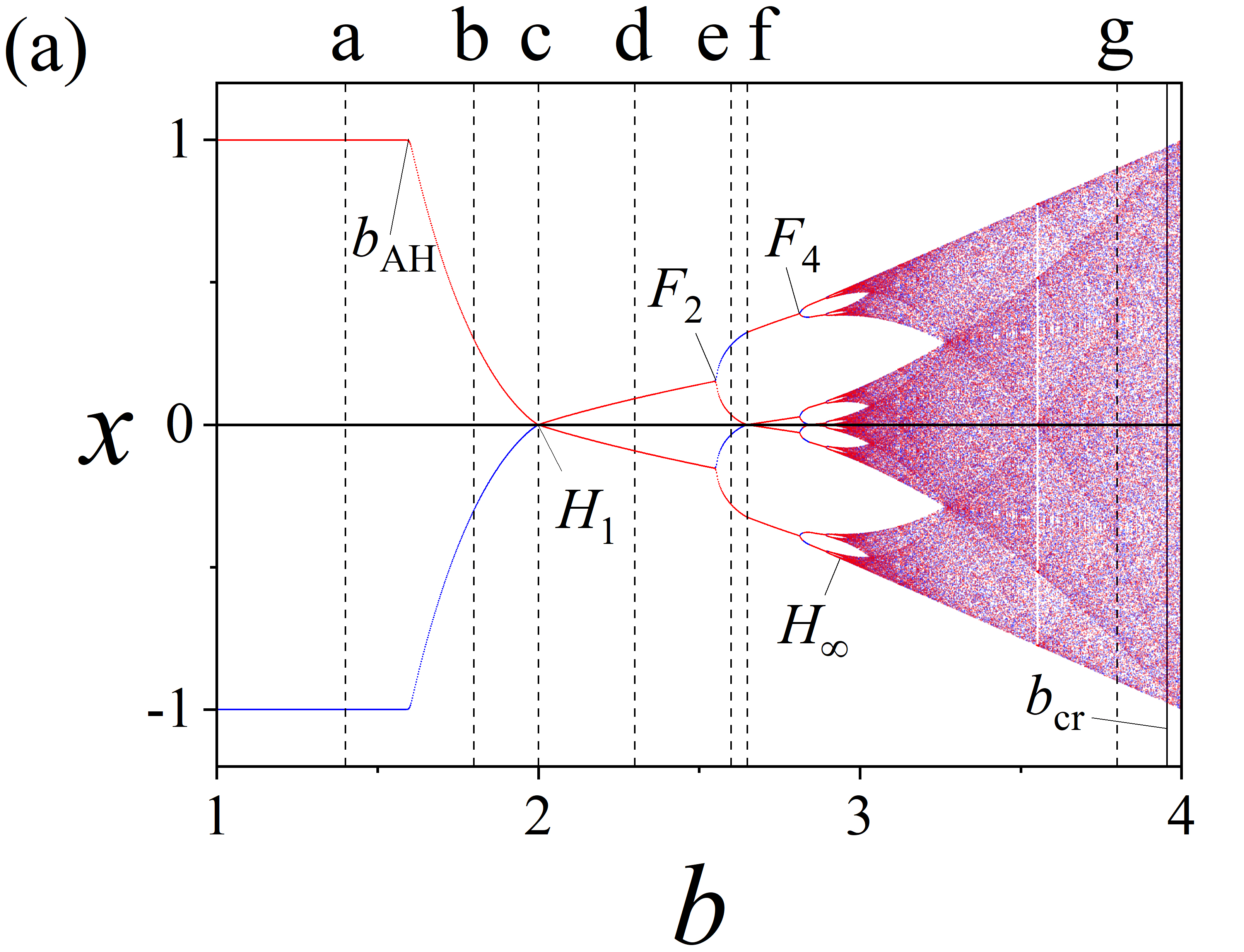}\quad  \includegraphics[width=0.47\textwidth]{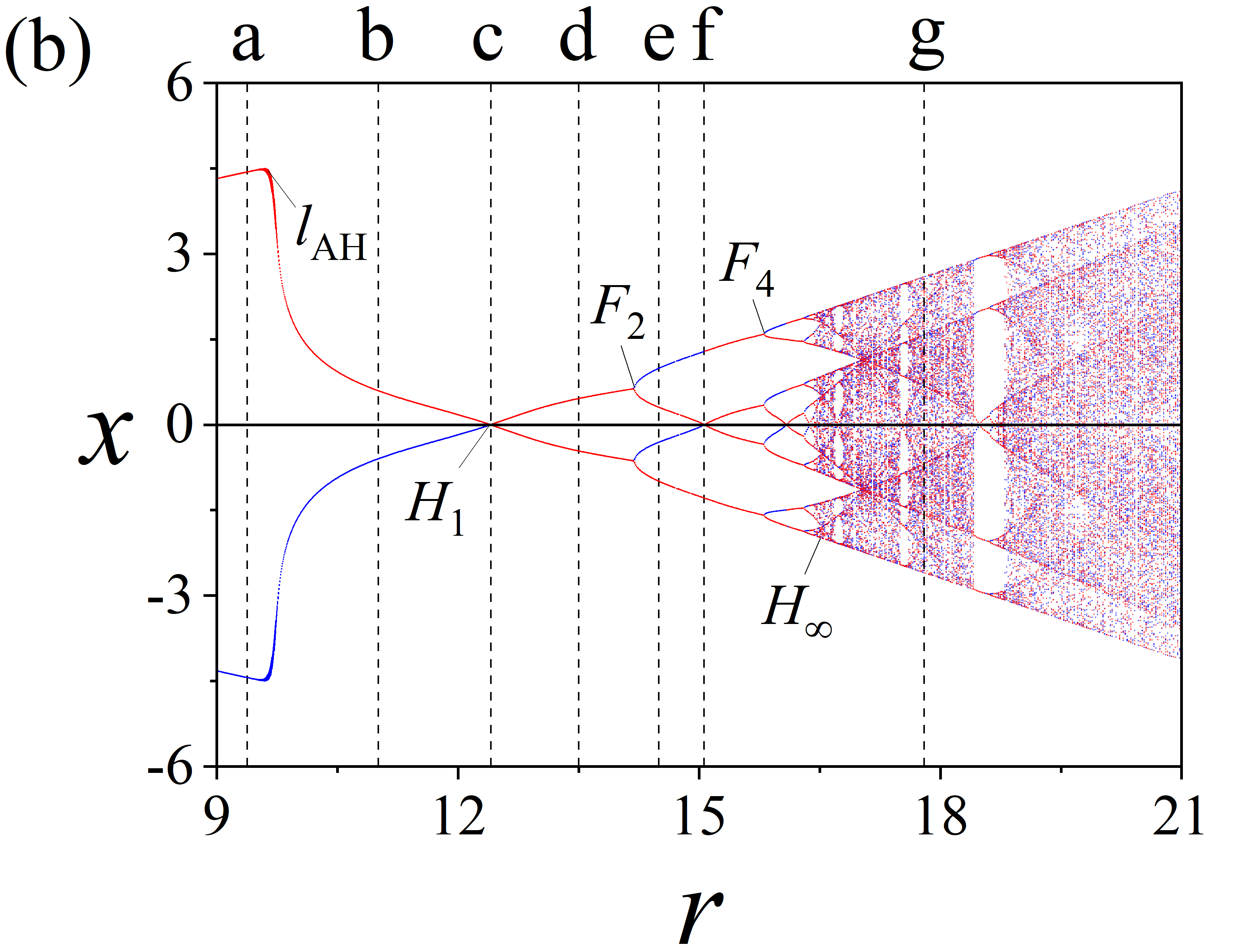}
    \caption{Extended 1-D bifurcation diagrams corresponding to the route $\text{a}-\text{f}$ in (a) PWL system~\eqref{map:master} and
    (b) LLZ system~\eqref{sys:Lyubimov}. 
The red and blue dots represent traces with initial conditions in the vicinity of $e_r$ and $e_l$, respectively, and indicate the coincidence or difference of the attracting sets. 
Coordinates of the route points and parameters are the same as in Figs.~\ref{fig:map_v>1}-\ref{fig:PWL_Lubimov_compare}.}    \label{fig:PWL_Lorenz_biff_gamma}
\end{figure}

\begin{figure}[ht]
\centering
\includegraphics[width=0.8\textwidth]{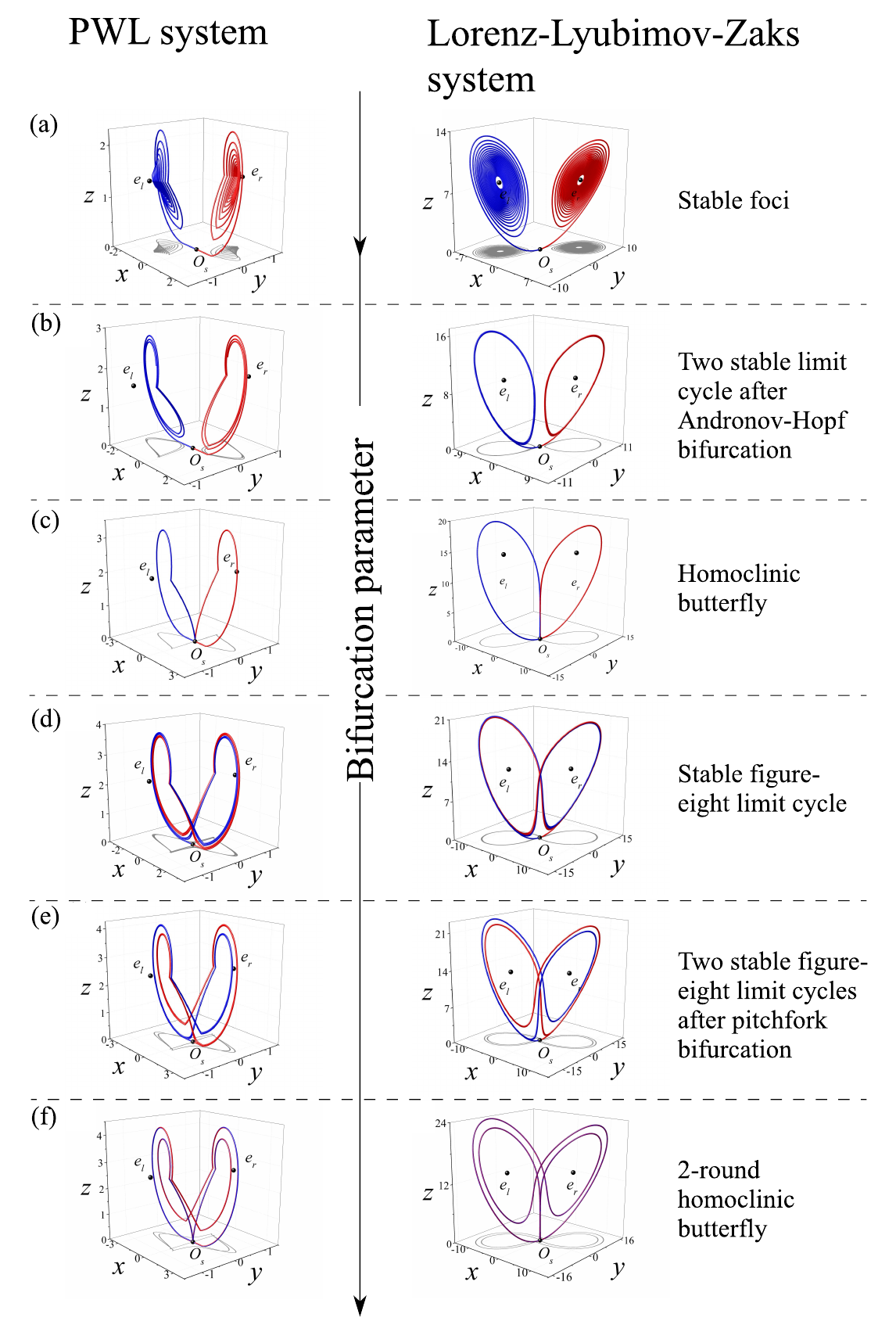}
\caption{Comparison the phase portraits of PWL system~\eqref{sys:PWL_Lorenz} and LLZ system~\eqref{sys:Lyubimov} along the route $\text{a}-\text{f}$ (see diagrams in Fig.~\ref{fig:PWL_Lorenz_biff_gamma}). The rows (a)--(f) correspond to the points $\text{a}-\text{f}$ in Fig.~\ref{fig:PWL_Lorenz_biff_gamma} and have the same parameters. Other parameters of PWL system~\eqref{sys:PWL_Lorenz}:
$\alpha = 2$, $\lambda = 0.294$, $\omega = 2$, $\delta = 0.588$.}
\label{fig:PWL_Lubimov_compare}
\end{figure}

The resulting expanded diagraMira is shown in Fig.~\ref{fig:PWL_Lorenz_biff_gamma}(a).

For LLZ system~\eqref{sys:Lyubimov}, we calculated the Poincar\'{e} return map along the flow trajectories at the following V-shape cross-section
\begin{equation*}
    \Pi: \quad z=\kappa^{-1}|\xi x+\eta y|,
\end{equation*}
where $\xi=(\sigma+s_2^0)(r+R)$, $\eta=\sigma(D-1)(r+R)$ and $\kappa=Q[R\sigma(D-1)-\sigma-s_2^0]$.
This cross-section passes through the saddle $O_s$, containing its strong stable eigenvector, and the foci $e_{l,r}$.

The parameters $(D,r)$ were changed according to route $\text{a}-\text{f}$ from Fig.~\ref{fig:2D_diags_negative}, i.e. $r=\frac{0.881}{D}$ [see Fig.~\ref{fig:PWL_Lorenz_biff_gamma}(b)].

Both diagrams in Fig.~\ref{fig:PWL_Lorenz_biff_gamma} show the similar bifurcation sets, including ``homoclinic-pitchfork'' cascade, the self-similarity phenomenon and periodicity windows. 
This further confirms that PWL system~\eqref{sys:PWL_Lorenz} closely reconstructs the main properties of quasi-strange Lorenz-type attractors.

\begin{figure}[h]
  \centering
\includegraphics[width=0.45\textwidth]{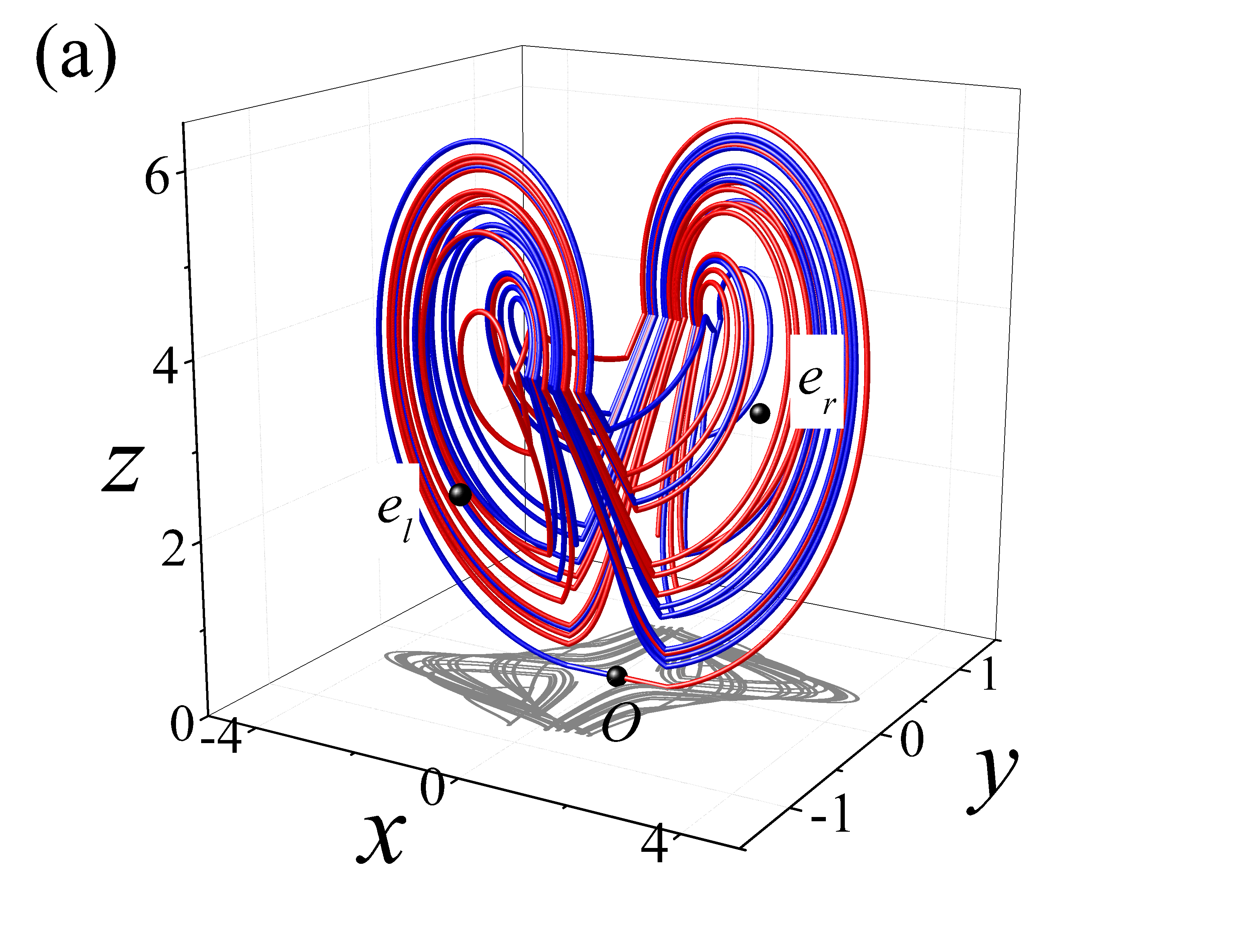}    \includegraphics[width=0.45\textwidth]{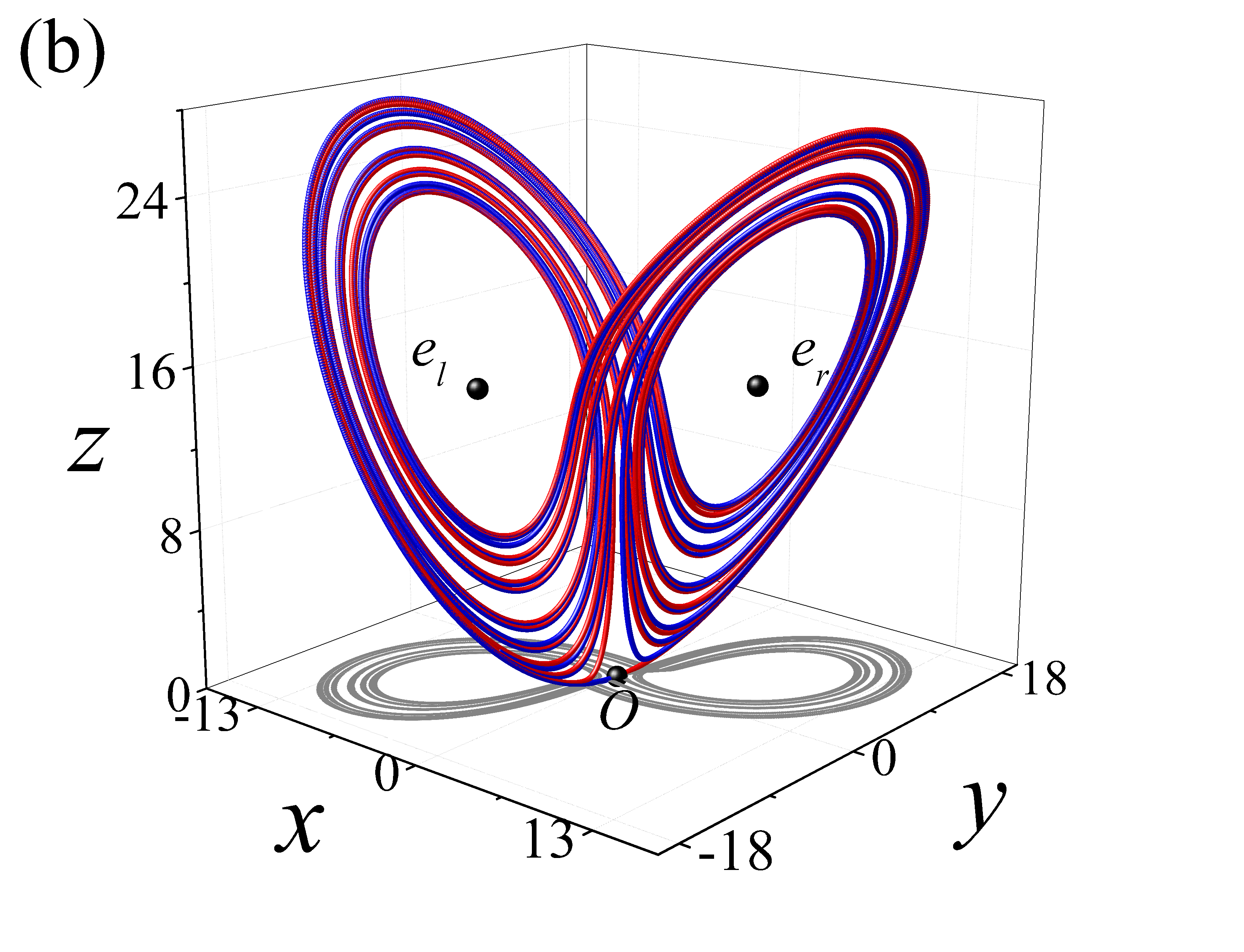}
    \caption{Phase portraits of the quasi-strange Lorenz-type attractor at the point $\text{h}$ in Fig.~\ref{fig:2D_diags_negative} for (a) PWL system~\eqref{sys:PWL_Lorenz} and (b) LLZ system~\eqref{sys:Lyubimov}.}
    \label{fig:Atractor_Lorenz_like}
\end{figure}

The phase portraits of such attractors for slice g from Fig.~\ref{fig:PWL_Lorenz_biff_gamma} are shown in Fig.~\ref{fig:Atractor_Lorenz_like}.

\textbf{Remark 1.}
Despite the similarity with the unimodal mapping, the systems~\eqref{sys:PWL_Lorenz}, \eqref{sys:Lyubimov} have a major difference.
Due to their odd symmetry, all periodic orbits are born either as a stable homoclinic butterfly $B^n$ via homoclinic bifurcation $H_n$, or as a symmetric pair $C^n_{1,2}$ as a result of the fold and pitchfork bifurcations. 
This concerns the stable birth of 3-periodic orbits $C^3_{1,2}$, which, unlike the smooth unimodal case, are born as a symmetric pair as a result of the fold bifurcation.

\textbf{Remark 2.}
A significant difference between the LLZ system~\eqref{sys:Lyubimov} and PWL system~\eqref{sys:PWL_Lorenz} is the possibility of sliding motions appearing in the latter, which radically changes the dynamics of the system.
In the system~\eqref{sys:PWL_Lorenz}, these motions appears for $b>b_{cr}$, and the PWL system no longer reproduces the dynamics of LLZ system~\eqref{sys:Lyubimov}.

%====================
% Thus, due to \eqref{eq:}, for $0<b<\frac{1}{\nu}e^{\frac{3\pi\lambda}{2\omega}}$, the system~\eqref{sys:PWL_Lorenz} has two stable foci, attracting the almost whole phase space (except 2-D stable manifolds of the saddle $O$).

% For $0<\gamma<1$

% When $\gamma = 1$, the point $x = 0$ is mapped into itself, which corresponds to the birth of the homoclinic butterfly $B^1$ in the system Eq.~\eqref{sys:PWL_Lorenz}.

% From the analysis of the mapping Eq.~\eqref{map:master}, we obtain the formulas of bifurcation curves.
% The curve equation for pitchfork bifurcation $F^2$ has the form
% \begin{equation*}
%     \left( \frac{1}{\nu \gamma} \right) ^ {\frac{1}{\nu - 1}} \left( \frac{1}{\nu} + 1 \right) = \gamma + 1,
% \end{equation*}
% and the equation of the homoclinic bifurcation curve $H^2$ corresponding to the birth of the homoclinic butterfly $B^2$ in the system Eq.~\eqref{sys:PWL_Lorenz} has the following form
% \begin{equation}
%     \nu = \log_{\gamma - 1} \frac{\gamma - 1}{\gamma}.
% \end{equation}

\section{Conclusion}
In this paper, we have analytically and numerically considered 3-D piecewise-linear system~\eqref{sys:PWL_Lorenz} of the Lorenz type.
We have shown that, depending on the value of the saddle index $\nu$, which is explicitly included in the system as a parameter, it qualitatively reproduces the key bifurcations and attractors of its two smooth originals - the Lorenz system with a singular hyperbolic attractor and the Lorenz-Lyubimov-Zaks system~\eqref{sys:Lyubimov} with a quasi-strange attractor (often called a ``contracting Lorenz attractor'').

Unlike smooth systems, for which the main analytical results are obtained in a small neighborhood of the homoclinic butterfly to the neutral saddle, a piecewise smooth system allowed us to analytically obtain a global partition of the parameter space, without resorting to asymptotic methods, assumptions and postulation of the necessary properties.

In particular, for case $\nu<1$ it is possible to obtain a heteroclinic bifurcation of codimension 1, leading to the birth of a singular hyperbolic attractor of the Lorenz type.

For case $\nu>1$, we showed the existence of a ``homoclinic-pitchfork''-cascade of codimension 1 bifurcations, the limit of which is the quasi-strange attractor.

In this paper, we consider bifurcations beyond this limit only numerically.
However, due to the Lorenz type of the factor map~\eqref{map:master}, the system~\eqref{sys:PWL_Lorenz} has the entire bifurcation structure of such mappings \cite{Malkin2003}.
It includes self embeddings similar to Mira's "box-within-a-box" bifurcation structure~\cite{Mira1987}, Sharkovsky type sequences~\cite{anuvsic2020topological} and period adding cascades~\cite{belykh1997comments,belykh20241d}.
These details are the subject of further consideration.

We hope that the approach to reconstructing nonlinear systems with complex dynamics presented here and in a number of other works~\cite{belykh2019lorenz,belykh2021sliding,belykh2023hidden} allows, firstly, to significantly improve and clarify the understanding of their dynamics, and, secondly, to reveal new effects that were unnoticed due to the cumbersomeness of the analysis methods or were inaccessible to them.

\section*{Acknowledgments}
Analytical results were supported by Russian Science Foundation (grant No. 24-71-00073), numerical results were supported by Ministry of Science and Education of Russian Federation (contract No. FSWR-2020-0036).

%\section{References}
% References cited in the text should be placed within square
% brackets and stated as [surname of author(s), year of
% publication], e.g., \cite{Golub89} and, with three
% or more authors, \cite{Hall97}. If the reference reads as part of
% the sentence, the square brackets enclose only the year of
% publication, e.g., ``According to Golub and Van Loan [1989], \ldots''

\noindent  Nikita V. Barabash - \url{https://orcid.org/0000-0001-7094-0776}

%\noindent Rajesh Babu - \url{https://orcid.org/0009-0006-0415-6880}

%\end{multicols}

\markboth{Authors' Names}{Paper Title}

\section*{References}

\begin{thebibliography}{99}

\markboth{Authors' Names}{Paper Title}

\bibitem[Lorenz(1963)]{lorenz1963deterministic}
Lorenz E.N. [1963]
`Deterministic Nonperiodic Flow,''
{\it Journal of Atmospheric Sciences} {\bf 20}, pp. 130--141,
\url{https://doi.org/10.1175/1520-0469(1963)020<0130:DNF>2.0.CO;2}

\bibitem[Guckenheimer(1976)]{guckenheimer1976strange}
Guckenheimer, J. [1976] ``A strange, strange attractor, The Hopf bifurcation and its applications,'' {\it Applied Mathematical Series} {\bf 19} (Springer-Vedlag), pp. 368--381.
\url{https://doi.org/10.1007/978-1-4612-6374-6_25}

\bibitem[Afraimovich {\it et al.}(1982)]{afraimovich1982attractive}
Afraimovich, V. S., Bykov, V. V. \& Shilnikov, L. P. [1982]
``Attractive nonrough limit sets of Lorenz-attractor type,'' 
{\it Trudy Moskovskoe Matematicheskoe Obshchestvo} {\bf 44}, pp. 150--212.

\bibitem[Afraimovich {\it et al.}(1977)]{afraimovich1977origin}
Afraimovich, V. S., Bykov, V. V. \& Shilnikov, L. P. [1977] ``On the origin and structure of the Lorenz attractor,'' {\it Akademiia Nauk SSSR Doklady} {\bf 234}, pp. 336--339.  

\bibitem[Williams(1979)]{williams1979structure}
Williams, R. F. [1979] ``The structure of Lorenz attractors,'' {\it Publications Mathematiques de l'IHES} {\bf 50}, pp. 73--99, \url{https://doi.org/10.1007/BF02684770}

\bibitem[Filippov(2013)]{filippov2013differential}
Filippov A.F. [2013]
{\it Differential equations with discontinuous righthand sides: control systems}, (Springer Science \& Business Media, USA).

\bibitem[Shilnikov(1980)]{shilnikov1980appendix}
Shilnikov L. [1980]
`Appendix to Russian edition of The Hopf Bifurcation and Its Applications,'
{\it The Hopf Bifurcation and Its Applications}, eds.~J. Marsden and M. McCraken, vol. 317.

\bibitem[Sparrow(1982)]{sparrow1982lorenz}
Sparrow, C., \textit{The lorenz equations}, Springer New York, 1982.

\bibitem[Petrovskaya(1980)]{petrovskaya1980homoclinic}
Petrovskaya, N. V., \& Yudovich, V. I. [1980] ``Homoclinic loops of the Saltzman-Lorenz system,'' 
{\it Methods of qualitative theory of differential equations} {\bf 20}, pp. 73--83.


\bibitem[Belykh(1984)]{belykh1984bifurcation}
Belykh, V. N. [1984]
``Bifurcation of separatrices of a saddle-point of the Lorenz system,'' 
{\it Differential equations} {\bf 20}, pp. 1184--1191.

\bibitem[Ovsyannikov(2016)]{ovsyannikov2016analytic}
Ovsyannikov, I. I. \& Turaev, D. V. [1963] ``Analytic proof of the existence of the Lorenz attractor in the extended Lorenz model,'' {\it Nonlinearity} {\bf 30}, pp. 115,
\url{http://doi.org/10.1088/1361-6544/30/1/115}

\bibitem[Tucker(1999)]{tucker1999lorenz}
Tucker, W. [1999] ``The Lorenz attractor exists,'' {\it Comptes Rendus de l'Académie des Sciences-Series I-Mathematics} {\bf 328}, pp. 1197--1202, \url{https://doi.org/10.1016/S0764-4442(99)80439-X}

\bibitem[Bykov~\& Shilnikov, 1989]{bykov1989boundaries}
Bykov, V. V., and Shilnikov, A. L., ``On the boundaries of the domain of existence of the Lorenz attractor,'' \textit{Methods of Qualitative Theory and Theory of Bifurcations. Gorky State University, Gorky}, 1989, pp. 151--159.

\bibitem[Barrio~\& Shilnikov, 2012]{barrio2012kneadings} Barrio, R., Shilnikov, A., and Shilnikov, L., Kneadings, symbolic dynamics and painting Lorenz chaos, \textit{International Journal of Bifurcation and Chaos}, 2012, {\bf 22}, pp. 1230016.

\bibitem[Creaser {\it et al.}(2017)]{creaser2017finding}
Creaser, J. L., Krauskopf, B., and Osinga, H. M., ``Finding first foliation tangencies in the Lorenz system,'' \textit{SIAM Journal on Applied Dynamical Systems}, 2017, {\bf 16}, pp. 2127--2164.

\bibitem[Arneodo {\it et al.}(1981)]{arneodo1981possible}
Arneodo, A., Coullet, P., \& Tresser, Ch. [1981] ``A possible new mechanism for the onset of turbulence,'' {\it Physics Letters A} {\bf 81}, pp. 197--201, \url{https://doi.org/10.1016/0375-9601(81)90239-5}

\bibitem[Lyubimov~\& Zaks, 1983]{lyubimov1983two}
Lyubimov, D. V. and Zaks, M. A., ``Two mechanisms of the transition to chaos in finite-dimensional models of convection,'' \textit{Physica D: Nonlinear Phenom.}, 1983, {\bf 9}, pp. 52--64.

\bibitem[Belykh {\it et al.}(2019)]{belykh2019lorenz}
Belykh, V. N., Barabash, N. V., \& Belykh, I. [2019] ``A Lorenz-type attractor in a piecewise-smooth system: Rigorous results,'' {\it Chaos} {\bf 29}, pp. 103108, \url{https://doi.org/10.1063/1.5115789}

\bibitem[Belykh {\it et al.}(2020)]{belykh2020bifurcations}
Belykh V.N., Barabash N.V., \& Belykh I. [2020]
``Bifurcations of chaotic attractors in a piecewise smooth Lorenz-type system,''
{\it Automation and Remote Control} {\bf 81}, pp. 1385--1393,
\url{https://doi.org/10.1134/S0005117920080020}

\bibitem[Belykh {\it et al.}(2021)]{belykh2021sliding}
Belykh V.N., Barabash N.V., \& Belykh I. [2021]
``Sliding homoclinic bifurcations in a Lorenz-type system: Analytic proofs,''
{\it Chaos} {\bf 31}, pp. 043117,
\url{https://doi.org/10.1063/5.0044731}

\bibitem[Belykh {\it et al.}(2023)]{belykh2023hidden}
Belykh V.N., Barabash N.V., \& Belykh I. [2023]
``The hidden complexity of a double-scroll attractor: Analytic proofs from a piecewise-smooth system,''
{\it Chaos} {\bf 33}, pp. 043119,
\url{https://doi.org/10.1063/5.0139064}

\bibitem[Rovella (1993)]{rovella1993dynamics}
Rovella A. [1993]
``The dynamics of perturbations of the contracting Lorenz attractor,''
{\it Bulletin/Brazilian Mathematical Society} {\bf 24}, pp. 233--259,
\url{https://doi.org/10.1007/BF01237679}

\bibitem[Morales {\it et al.}(2006)]{morales2006contracting}
Morales, C. A., Pacifico, M. J., \& Martin, B. S. [2006] ``Contracting Lorenz attractors through resonant double homoclinic loops,'' {\it SIAM Journal on Mathematical Analysis} {\bf 38}, pp. 309--332, \url{https://doi.org/10.1137/S0036141004443907}

\bibitem[Anu{\v{s}}i{\'c}(2020)]{anuvsic2020topological}
Anu{\v{s}}i{\'c}, A., Bruin, H., \& {\v{C}}in{\v{c}} J. [2020] ``Topological properties of Lorenz maps derived from unimodal maps,'' {\it Journal of Difference Equations and Applications} {\bf 26}, pp. 1174--1191, \url{https://doi.org/10.1080/10236198.2020.1760260}

\bibitem[Araujo(2021)]{araujo2021statistical}
Araujo V.[2021] ``On the statistical stability of families of attracting sets and the contracting Lorenz attractor,''
{\it Journal of Statistical Physics} {\bf 182}, pp. 53, 
\url{https://doi.org/10.1007/s10955-021-02729-x}

\bibitem[Kazakov(2021)]{kazakov2021bifurcations}
Kazakov A. [2021] ``On bifurcations of Lorenz attractors in the Lyubimov--Zaks model,''
{\it Chaos} {\bf 31}, pp. 093118, \url{https://doi.org/10.1063/5.0058585}

\bibitem[Li~\& Malkin, 2003]{Malkin2003}
Li, M. -C., and Malkin, M., ``Smooth symmetric and Lorenz models for unimodal maps,'' \textit{Int. J. Bifurcation Chaos}, 2003, {\bf 13}, pp. 3353--3371.

\bibitem[Mira(1987)]{Mira1987}
Mira C. [1987]
{\it Chaotic dynamics: From the one-dimensional endomorphism to the two-dimensional diffeomorphism}, (World Scientific).

\bibitem[Belykh~\& Mosekilde(1997)]{belykh1997comments}
Belykh V.N., \& Mosekilde E. [1997]
``Comments on the Theory of Unimodal Maps,''
{\it Open Systems \& Information Dynamics} {\bf 4}, pp. 379--392,
\url{https://doi.org/10.1023/A:1009664631952}

\bibitem[Belykh {\it et al.}(2024)]{belykh20241d}
Belykh V.N., Barabash N.V., \& Grechko D.A.[2024]
``From 1D endomorphism to multidimensional H{\'e}non map: Persistence of bifurcation structure,''
{\it International Journal of Bifurcation and Chaos} {\bf 34}, pp. 2450026,
\url{https://doi.org/10.1142/S0218127424500263}




% \bibitemira[Zhu~\& Leung, 1999]{EKF_DrLeung3}
% Zhu, Z., Leung, H. \& Ding, Z. [1999]
% ``Optimiraal synchronization of chaotic systems in noise,''
% {\it IEEE Trans. Circ. Syst.-I$:$ Fund. Th. Appl.} {\bf 46}, 1320--1329,
% \url{https://doi.org/10.1109/81.802822}.

% \bibitem[Haller {\it et al.}(1997)]{Hall97}
% Haller, B., Streiff, M., Fleisch, U. \& Zimmermann, R. [1997]
% ``Hardware implementation of a systolic \hbox{antenna} array signal processor based on CORDIC arithmetic,''
% {\it Proc. Int. Conf. Acoust. Speech \hbox{Signal} Proc.} {\bf 5}, pp. 4141--4144, \url{https://doi.org/10.1109/ICASSP.1997.604858}

% \bibitem[Golub \& Van Loan(1989)]{Golub89}
% Golub, G. H. \& Van Loan, C. F. [1989]
% {\it Matrix Comptations}, 2nd Ed. (Johns Hopkins University Press, USA).

% \bibitem[Haller {\it et al.}(1997)]{Hall97}
% Haller, B., Streiff, M., Fleisch, U. \& Zimmermann, R. [1997]
% ``Hardware implementation of a systolic \hbox{antenna} array signal processor based on CORDIC arithmetic,''
% {\it Proc. Int. Conf. Acoust. Speech \hbox{Signal} Proc.} {\bf 5}, pp. 4141--4144, \url{https://doi.org/10.1109/ICASSP.1997.604858}.

% \bibitem[Lie \& Wang(2000)]{Lie00}
% Lie, D. Y. C. \& Wang, K. L. [2000]
% ``Si/SiGe heterostructures for Si-based nanoelectronics,''
% {\it Handbook of} \hbox{\it Advanced} {\it Electronic and Photonic Devices and Materials}, ed.~Nalwa, H. S. (\hbox{Academic} Press, San Diego), Chapter~1, pp.~1--69,
% \url{https://doi.org/10.1016/B978-012513745-4/50021-4}.

% \bibitem[Lie \& Wang(2001)]{Lie01}
% Lie, D. Y. C. \& Wang, K. L. [2001]
% ``SiGe/Si processing,'' {\it Semiconductors and Semimetals} {\bf 73},
% eds.~Willardson, R. \& Weber, E. (Academic Press, San Diego), Chapter 4, pp.~151--197,
% \url{https://doi.org/10.1016/S0080-8784(01)80214-3}.

% \bibitem[Parssinen {\it et al.}(1999)]{Par99}
% P\"arssinen, A., Jussila, J., Ryyn\"anen, J., Sumanen, L. \& Halonen, K.  A. I. [1999]
% ``A 2-GHz wide-band direct conversion receiver for WCDMA applications,''
% {\it IEEE J. Solid-State Circuits} {\bf 34}, 1893,
% \url{https://doi.org/10.1109/4.808914}.

\end{thebibliography}
\end{document}